\documentclass{article}
\usepackage{amsmath,amsxtra,amssymb,latexsym,amsfonts,amscd}
\usepackage{multicol,color}
\usepackage{float}
\usepackage{soul}
\usepackage{graphicx}
\definecolor{labelkey}{rgb}{0,0.08,0.45}
\definecolor{refkey}{rgb}{0,0.6,0.0}
\definecolor{Brown}{rgb}{0.45,0.0,0.05}
\definecolor{lime}{rgb}{0.00,0.8,0.0}
\definecolor{lblue}{rgb}{0.5,0.5,0.99}
\usepackage{pstricks-add}
\usepackage{amssymb}
\usepackage{theorem}
\usepackage{fancyhdr}
\usepackage{tikz}
\usepackage{enumerate}
\usepackage[margin=1.1in]{geometry}
\usetikzlibrary{arrows}
\usepackage{hyperref}
\definecolor{labelkey}{rgb}{0.6,0.6,0.6}
\definecolor{refkey}{rgb}{0,0.6,0.0}
\def\disp{\displaystyle}
\def\e{\varepsilon}
\def\ph{\varphi}
\def\ve{\varepsilon}
\def\dd{\delta}

\def\lm{\lambda}
\def\O{\Omega}
\def\Tilde{\widetilde}

\def\tx{\widetilde{x}}
\def\tu{\widetilde{u}}
\def\ta{\widetilde{a}}
\def\tb{\widetilde{b}}

\def\({\left(}
\def\){\right)}
\def\[{\left[}
\def\]{\right]}
\def\n{\left \|}
\def\en{\right \|}
\def\nn{\left \{ }
\def\hnn{\right \}}

\def\oa{\bar a}
\def\ob{\bar b}

\def\ox{\bar{x}}
\def\oy{\bar{y}}
\def\oz{\bar{z}}
\def\ov{\bar{v}}
\def\ou{\bar{u}}

\def\oA{\bar A}
\def\oB{\bar B}
\def\oX{\bar X}
\def\oal{\bar\alpha}
\def\obe{\bar\beta}

\def\gph{\hbox{}}

\def\gg{\gamma}
\def\dn{\downarrow}

\def\sr{\Longrightarrow }
\def\tto{\rightrightarrows}

\def\Limsup{\mathop{{\rm Lim}\,{\rm sup}}}

\def\hat{\widehat}
\def\Hat{\widehat}

\def\Tilde{\widetilde}
\def\Bar{\overline}
\def\la{\left\langle}
\def\ra{\right\rangle}
\def\ve{\varepsilon}

\def\h{\hfill\Box}
\def\R{\mathbb{R}}
\def\N{\mathbb{N}}

\def\gra{\bigtriangledown}
\def\co{\mbox{\rm co}\,}
\def\gph{\mbox{\rm gph}\,}
\def\rep{\mbox{\rm rep}\,}
\def\epi{\mbox{\rm epi}\,}

\def\dom{\mbox{\rm dom}\,}

\def\var{\mbox{\rm var}\,}
\def\dn{\downarrow}
\def\O{\Omega}

\def\vph{\varphi}
\def\emp{\emptyset}

\def\oR{\Bar{\R}}
\def\lm{\lambda}

\def\gg{\gamma}

\def\dd{\delta}

\def\al{\alpha}

\def\t{\tau}

\def\N{I\!\!N}
\def\th{\theta}

\newtheorem{theorem}{Theorem}[section]
\newtheorem{lemma}[theorem]{Lemma}

\newtheorem{proposition}[theorem]{Proposition}
\newtheorem{definition}[theorem]{Definition}
\theoremstyle{plain}{\theorembodyfont{\rmfamily}
}
\theoremstyle{plain}{\theorembodyfont{\rmfamily}
}
\theoremstyle{plain}{\theorembodyfont{\rmfamily}
}
\theoremstyle{plain}{\theorembodyfont{\rmfamily}
}

\theoremstyle{plain}{\theorembodyfont{\rmfamily}
\newtheorem{remark}[theorem]{Remark}}

\def\eq{\begin{equation}}
\def\eeq{\end{equation}}
\begin{document}
\begin{center}
{\bf OPTIMIZATION AND DISCRETE APPROXIMATION OF SWEEPING PROCESSES WITH
CONTROLLED MOVING SETS AND PERTURBATIONS}
\\[3ex]
TAN H. CAO\footnote{Department of Applied Mathematics and Statistics, State
University of New York--Korea, Yeonsu-Gu, Incheon, Republic of Korea
(tan.cao@stonybrook.edu). Research of this author is partially supported by IFA (Individual Faculty Accounts) funding.}
\quad GIOVANNI COLOMBO\footnote{Dipartimento di Matematica ``Tullio Levi-Civita",
Universit$\grave{\textrm{a}}$ di Padova, via Trieste 63, 35121 Padua, Italy
(colombo@math.unipd.it). Research of this author is partially supported by the University of Padova grant SID 2018 ``Controllability, stabilizability and infimun gaps for control systems," BIRD 187147, and is affiliated to Istituto Nazionale di Alta Matematica (GNAMPA).}
\quad BORIS S. MORDUKHOVICH\footnote{Department of Mathematics, Wayne State
University, Detroit, Michigan 48202, USA (boris@math.wayne.edu). Research of
this author was partially supported by the USA National Science Foundation under
grants DMS-1007132 and DMS-1512846, by the USA Air Force Office of Scientific
Research grant \#15RT0462, and by the Australian Research Council under
Discovery Project DP-190100555.}
\quad DAO NGUYEN\footnote{Department of
Mathematics, Wayne State University, Detroit, Michigan 48202, USA
(dao.nguyen2@wayne.edu). Research of this author was partially supported by
the USA National Science Foundation under grant DMS-1512846 and by the USA Air
Force Office of Scientific Research grant \#15RT0462.}
\end{center}
\small{\bf Abstract.}
This paper addresses a new class of optimal control problems for perturbed
sweeping processes with measurable controls in additive perturbations of the
dynamics and smooth controls in polyhedral moving sets. We develop a
constructive discrete approximation procedure that allows us to strongly
approximate any feasible trajectory of the controlled sweeping process by
feasible discrete trajectories and also establish a $W^{1,2}$-strong convergence
of optimal trajectories for discretized control problems to a given local
minimizer of the original continuous-time sweeping control problem of the Bolza
type. Employing advanced tools of first-order and second-order variational
analysis and generalized differentiation, we derive necessary optimality
conditions for discrete optimal solutions under fairly general assumptions
formulated entirely in terms of the given data. The obtained results give us
efficient suboptimality (``almost optimality") conditions for the original
sweeping control problem that are illustrated
by a nontrivial numerical example.\\[0.02in]
{\em Key words.} Optimal control, Sweeping process, Discrete approximation,
Variational analysis, Generalized differentiation, Necessary optimality
conditions\\[0.02in]
{\em AMS Subject Classifications}: 49J52; 49J53; 49K24; 49M25;
90C30\vspace{-0.2in}

\section{Problem Formulation and Initial Discussions}\label{intro}
\setcounter{equation}{0}\vspace*{-0.1in}

This paper is devoted to the study of {\em optimal control} problems for {\em
sweeping processes} with {\em controlled perturbations} and {\em controlled
moving sets}. The basic {\em uncontrolled} sweeping process was introduced by
Moreau in the 1970s as the dissipative differential inclusion
\begin{equation}\label{SP}
\dot x(t)\in-N\big(x(t);C(t)\big)\;\mbox{ a.e. }\;t\in[0,T]\;\mbox{ with
}\;x(0):=x_0\in C(0)
\end{equation}
describing the motion of a particle that belongs to a continuously moving set
$C(t)$, where the normal cone $N$ in \eqref{SP} is understood in the sense of
convex analysis
\begin{equation}\label{NC}
N(x;C)=N_C(\ox):=\big\{v\in\R^n\big|\;\la v,y-x\ra\le 0,\;y\in
C\big\}\;\textrm{if}\;x\in C\;\textrm{ and }\;N(x;C):=\emp\textrm{ if }x\notin
C.
\end{equation}
The sweeping inclusion \eqref{SP} tells us that, depending on the motion of the
set, the particle stays where it is in the case when it does not hit the set;
otherwise, it is swept towards the interior of the set. We refer the reader to
\cite{mor_frict} and to the subsequent work in, e.g.,
\cite{aht,BrS,CT,cpt,jv,KP,KMM,mv,ve} with the bibliographies therein for
further developments and applications. The original motivation for Moreau came
from applications to elastoplasticity, but later on the sweeping process and its
modifications have been well recognized for many applications to other problems
in mechanics, hysteresis, ferromagnetism, electric circuits, phase transitions,
traffic equilibria, social and economic modelings, etc.; see, e.g., the
references above among numerous publications.\vspace*{-0.05in}

Since the Cauchy problem in \eqref{SP} has a {\em unique} solution
\cite{mor_frict}, it does not make any sense to formulate optimization problems
for the basic Moreau sweeping process. This is a striking difference between the
discontinuous differential inclusion \eqref{SP} and the ones $\dot x(t)\in
F(x(t))$ described by {\em Lipschitzian} set-valued mappings/multifunctions
$F\colon\R^n\tto\R^n$ for which optimal control theory has been well developed;
see, e.g., the books \cite{c,m-book,v} for various methods and results on
necessary optimality conditions.\vspace*{-0.05in}

It seems that optimal control problems for sweeping differential inclusions were
first formulated and studied in the case of control actions entering {\em
additive perturbations} \cite{et} for which existence and relaxation results,
while not optimality conditions, were obtained; see \cite{cmr,Tol1,Tol} for
subsequent developments in this direction. To the best of our knowledge, the
theory of necessary optimality conditions for sweeping processes has been
started with \cite{chhm1}, where a new class of dynamic optimization problems
with {\em controlled moving sets} $C(t)=C(u(t))$ in \eqref{SP} was first
formulated with deriving necessary optimality conditions in the case when $C(u)$
is a half-space. Soon after that, necessary optimality conditions were obtained
for another class of sweeping process without controlled in either moving sets
or perturbations, but in a coupling linear ODE. Further necessary optimality
conditions and their applications for all the three types of controlled sweeping
processes were
developed in \cite{ao,ac,cm1,cm2,cm3,cm4,chhm2,chhm3,cmn1,cmn,pfs,hm}.

This paper concerns the following class of optimal control problems of the {\em
generalized Bolza type} for the perturbed version of the sweeping process in
\eqref{SP}. Given an extended real-valued terminal cost function
$\vph\colon\R^n\to\overline{\R}:=\(-\infty,\infty\]$ and a running cost function
$\ell\colon\[0,T\]\times\R^{2(n+nm+m)+d}\to\overline{\R}$, our basic problem
$(P)$ is defined by:
\begin{equation}\label{minimize}
\textrm{minimize}\hspace{0.1in}J[x,a,b,u]:=\vph\big(x(T)\big)+\int_0^T\ell\big(t
,x(t),a(t),b(t),u(t),\dot{x}(t),\dot{a}(t),\dot{b}(t)\big)dt
\end{equation}
over control actions $a(\cdot)=(a_1(\cdot),\ldots,a_m(\cdot))\in
W^{1,2}([0,T];\R^{mn})$ and $b(\cdot)=(b_1(\cdot),\ldots,b_m(\cdot))\in
W^{1,2}([0,T];{\R^m})$ entering the moving set $C(t)$ and measurable controls
$u(\cdot)\in L^2([0,T];\R^d)$ entering additive perturbations that generate the
corresponding trajectories $x(\cdot)\in W^{1,2}([0,T];\R^n)$ of the sweeping
differential inclusion
\begin{equation}\label{Problem}
\left\{\begin{array}{ll}
\dot{x}(t)\in-N\big(x(t);C(t)\big)+g\big(x(t),u(t)\big)\;\textrm{ a.e.
}\;t\in[0,T],\\
x(0):=x_0\in C(0),\;u(t)\in U\;\textrm{ a.e. }\;t\in[0,T],
\end{array}\right.
\end{equation}
where the moving set is given in the polyhedral form as
\begin{equation}\label{C}
C(t):=\big\{x\in\R^n\big|\;\la a_i(t),x\ra\le b_i(t),\;i=1,\ldots,m\big\},
\end{equation}
and where the initial point $x_0\in\R^n$ and the final time $T>0$ are fixed. All
such quadruples $(x(\cdot),a(\cdot),b(\cdot),u(\cdot))$ for which the running
cost $\ell(\cdot)$ is integrable are {\em feasible solutions} to problem $(P)$.
\vspace*{-0.05in}

In addition to the above dynamical system \eqref{Problem} with the {\em
pointwise/hard constraints} on the controls $u(\cdot)$ in perturbations, we
impose the pointwise constraints on the controls $a_i(\cdot)$ in the moving set:
\begin{equation}\label{h1.5}
\|a_i(t)\|=1\;\textrm{ for all }\;t\in [0,T]\;\mbox{ and }\;i=1,\ldots,m.
\end{equation}
Furthermore, problem $(P)$ also contains the implicit pointwise {\em mixed
state-control constraints}
\begin{equation}\label{e:8}
\la a_i(t),x(t)\ra\le b_i(t)\;\textrm{ for all }\;t\in[0,T]\;\textrm{ and
}\;i=1,\ldots,m,
\end{equation}
which are due to construction \eqref{NC} of the normal cone in
\eqref{Problem}.\vspace*{-0.05in}

Our approach to the dynamic optimization problem $(P)$ is based on the {\em
method of discrete approximation}, which was developed in \cite{m95,m-book} for
optimization of Lipschitzian differential inclusions and then was significantly
modified in \cite{cm1,cm2,cm3,chhm1,chhm2,chhm3,cmn,hm} to handle various
optimal control problems for sweeping processes. There are {\em four major
steps} in the realization of this approach to the study of continuous-time
systems:\\[1ex]
{\bf(i)} Firstly, we construct a {\em well-posed} discrete approximation of the
sweeping control system from \eqref{SP}, \eqref{C} in such a way that {\em any
feasible solution} to the continuous-time sweeping inclusion can be {\em
appropriately approximated} by feasible solutions to the discretized sweeping
control systems. This step may be also considered from the {\em numerical
viewpoint} as a finite-dimensional approximation of the discontinuous
constrained differential inclusion.\\[1ex]
{\bf(ii)} The second approximation step is to construct, with the usage if (i),
a sequence of discrete-time optimal control problems $(P_k)$,
$k\in\N:=\{1,2,\ldots\}$, for discretized sweeping inclusions such that the
approximating problems admit optimal solutions whose continuous-time extensions
{\em strongly converge} as $k\to\infty$ in the requited topology to a {\em
chosen local minimizer} of the original sweeping control problem $(P)$.\\[1ex]
{\bf(iii)} The next step is to derive {\em necessary conditions} that hold for
optimal solutions of each discrete-time problem $(P_k)$, which can be reduced to
a finite-dimensional format of {\em mathematical programming} with increasingly
many {\em geometric constraints} of the graphical type. To deal with such
problems, we employ appropriate tools of {\em first-order} and {\em second-order
variational analysis} and {\em generalized differentiation}. Due to (ii), the
obtained results can be viewed as constructive {\em suboptimality} (almost
optimality) conditions for $(P)$ that practically provide, for large $k\in\N$,
about the same amount of information as the exact optimality conditions for
local minimizers of $(P)$.\\[1ex]
{\bf(iv)} The last step is highly challenging mathematically while being of
undoubted importance. It furnishes the limiting procedure to pass from the
necessary conditions for the optimal solutions of the discrete-time problems
$(P_k)$ obtained in (iii) to the {\em exact necessary optimality conditions} for
the designated {\em local minimizer} of the original sweeping control problem
$(P)$. This step strongly involves advanced calculus and computation results of
variational analysis and generalized differentiation, especially of the {\em
second order}.
\vspace*{-0.02in}

In this paper we comprehensively resolve the issues listed in steps (i)--(iii)
for the general sweeping control problem $(P)$ formulated in
\eqref{minimize}--\eqref{e:8} (which is certainly of its independent interest
and own importance), while step (iv) is furnished in our forthcoming paper
\cite{ccmn}. Note that some particular cases of problem $(P)$ were investigated
by discrete approximation techniques in the papers \cite{cm1,cm3,chhm3,cmn}
mentioned above, but the general setting of our consideration is significantly
more complicated and thus requires careful elaborations, which are provided in
this paper and subsequently in \cite{ccmn}. \vspace*{-0.05in}

The rest of the paper is organized as follows. In Section~\ref{prel} we
formulate the {\em standing assumptions} on the given data of $(P)$ and present
preliminary results on the well-posedness of the controlled sweeping process
under consideration. Section~\ref{exist} establishes the {\em existence of
optimal solutions} to $(P)$ and discusses its relaxation stability. In
Section~\ref{disc} we construct a discrete approximation of the sweeping control
system in \eqref{SP}, \eqref{C} that allows us to {\em strongly approximate} any
{\em feasible solution} to it by feasible solutions to its discrete
counterparts. Section~\ref{disc1} develops the discrete approximation procedure
at the {\em level of optimality} while leading us to the strong convergence of
optimal solutions for the discrete-time problems to the prescribed {\em local
minimizer} of $(P)$. In Section~\ref{tools} we first review the tools of {\em
generalized differentiation} needed for our variational analysis and then obtain
{\em second-order
calculation formulas} that are crucial for deriving necessary optimality
conditions. Such conditions are obtained in Section~\ref{nec} for the
constructed discrete approximation problems, Finally, we illustrate in
Section~\ref{exam} by a nontrivial example the efficiency of the obtained
optimality conditions to solve sweeping control problems. Throughout the paper
we use standard notation of variational analysis and control theory; see, e.g.,
\cite{m18,rw,v}.\vspace*{-0.2in}

\section{Standing Assumptions and Preliminaries}\label{prel}
\setcounter{equation}{0}\vspace*{-0.1in}

In this section we present some results on well-posedness of the sweeping
differential inclusions in the aforementioned classes of feasible controls and
formulate the standing assumptions on problem $(P)$ that allow us to establish
further the main achievements of the paper.\vspace*{-0.05in}

Denoting by $d(x;\O)$ the distance between a given point $x\in\R^n$ and an
nonempty set $\O\subset\R^n$, observe first that the conventional assumption on
the moving set $C(t)$ ensuring the existence of absolutely continuous solutions
to the sweeping differential inclusion \eqref{Problem} is formulated as follows:
\begin{equation}\label{H}
|d\big(x;C(t)\big)-d\big(x;C(s)\big)|\le|v(t)-v(s)|\;\mbox{ for all
}\;t,s\in[0,T],
\end{equation}
where $v\colon[0,T]\to\R$ is an absolutely continuous function; see
\cite{CT,KMM} and the references therein. However, assumption \eqref{H} is
rather restrictive and may fail for polyhedral moving sets $C(t)$ as in
\eqref{C}, even in the case of half-spaces. An improvement of \eqref{H} ensuring
the existence of absolutely continuous solutions to \eqref{Problem} was obtained
in \cite{chhm2} with the verification of the imposed assumption in the case of
half-spaces $C(t)$ in \cite{chhm2} and then for general convex polyhedral sets
\eqref{C} in \cite{chhm3} under the linear independence constraint qualification
(LICQ) meaning that the vectors $\{a_i(t)\}$ are linearly independent for all
$t\in[0,T]$ along the active constraints. Following the approach of Tolstonogov
\cite{Tol1}, we derive below an advanced result on the existence and uniqueness
of $W^{1,2}$ solutions to \eqref{Problem} with the polyhedral moving sets
\eqref{C} generated by $W^{1,2}$ controls $(a_i(t),b_i(t))$ and measurable
controls $u(t)$
under a major assumption that is significantly weaker than LICQ. This result
justifies the well-posedness of the sweeping dynamical systems under
consideration, which is required for the subsequent study of the optimal control
problem $(P)$.\vspace*{-0.05in}

Now we formulate the {\em standing assumptions} of this paper that include those
ensuring the existence of the aforementioned solutions to the sweeping system
\eqref{Problem} and \eqref{C}.

{\bf(H1)} The control set $U$ from \eqref{Problem} is closed and bounded in
$\R^d$.

{\bf(H2)} The derivatives $(\dot a_i(t),\dot b_i(t))$ are uniformly bounded for
all $i=1,\ldots,m$ and a.e.\ $t\in[0,T]$ with the fixed initial points
$a_0:=(a_1(0),\ldots,a_m(0))$ and $b_0:=(b_1(0),\ldots,b_m(0))$.

{\bf(H3)} The perturbation mapping $g\colon\R^n\times\R^d\to\R^n$ is uniformly
Lipschitz continuous with respect to both variables $x$ and $u\in U$, i.e.,
there exists $L>0$ for which
\begin{equation}\label{e:gL}
\left\|g(x_1,u_1)-g(x_2,u_2)\right\|\le
L\left(\left\|x_1-x_2\right\|+\left\|u_1-u_2\right\|\right)\;\;\mbox{for
all}\;\;(x_1,u_1)\;\mbox{ and }\;(x_2,u_2)\in\R^n\times U.
\end{equation}
Furthermore, $g$ satisfies the sublinear growth condition
\begin{equation*}
\|g(x,u)\|\le M\(1+\|x\|\)\;\mbox{ for all }\;u\in U\;\mbox{ with some }\;M>0.
\end{equation*}

{\bf(H4)} There are functions $v_i(\cdot)\in W^{1,2}([0,T];\R)$ as
$i=1,\ldots,m$ such that
\begin{equation}\label{v-i}
\n a_i(t)-a_i(s)\en\le|v_i(t)-v_i(s)|\;\mbox{ for all }\;s,t\in[0,T]\;\mbox{ and
}\;i=1,\ldots,m.
\end{equation}
In addition, there exists a continuous function $\vartheta\colon[0,T]\to\R$ for
which $\disp\sup_{t\in[0,T]}\vartheta(t)<0$ and
\begin{equation}\label{Slater}
C^0(t):=\big\{x\in\R^n\big|\;\la
a_i(t),x\ra-b_i(t)<\vartheta(t),\;i=1,\ldots,m\big\}\ne\emp\;\mbox{ for all
}\;t\in[0,T].
\end{equation}

{\bf(H5)} The terminal cost $\vph\colon\R^n\to\oR$ is lower semicontinuous
(l.s.c.) while the running cost/integrand
$\ell\colon\R^{2(n+nm+m)+d}\to\bar{\R}$ is bounded from below and l.s.c.\ around
a given feasible solution to $(P)$ for a.e.\ $t\in[0,T]$. We also assume that
$\ell$ is a.e.\ continuous in $t$ and is uniformly majorized by a summable
function on $[0,T]$.

Before presenting the aforementioned well-posedness (existence and uniqueness)
theorem for the sweeping process in \eqref{Problem} and \eqref{C}, we discuss
the imposed condition \eqref{Slater} in (H4). Recall that the {\em positive
linear independence constraint qualification} (PLICQ) condition holds at $x\in
C(t)$ if
\begin{equation}\label{PLICQ}
\left[\sum_{i\in I(x,a(t),b(t))}\al_i a_i(t)=0,\;\al_i\ge
0\right]\Longrightarrow\big[\al_i(t)=0\;\;\mbox{for all}\;\;i\in
I\big(x,a(t),b(t)\big)\big],
\end{equation}
where the set of {\em active constraint indices} for \eqref{C} is defined by
\begin{equation}\label{e:AI}
I\big(x,a(t),b(t)\big):=\big\{i\in\{1,\ldots,m\}\big|\;\la
a_i(t),x\ra=b_i(t)\big\},\quad t\in[0,T].
\end{equation}
The essentially more restrictive {\em linear independence constraint
qualification} (LICQ) condition at $x\in C(t)$ used in \cite{chhm3} reads as
\eqref{PLICQ} with the replacement of $\al_i\ge 0$ by $\al_i\in\R$ therein.

It is easy to see the Slater-type condition \eqref{Slater} reduces to PLICQ  if
the polyhedron \eqref{C} does not depend on $t$, which is the case considered in
\cite{cmn}. In the general nonautonomous case, \eqref{Slater} may be stronger
than PLICQ \eqref{PLICQ} while being always
weaker than its LICQ counterpart. Note also that in our setting, \eqref{PLICQ}
corresponds to the {\em Mangasarian-Fromovitz constraint qualification}, which
is classical in nonlinear programming. Furthermore, imposing PLICQ at $x\in
C(t)$ is equivalent to the so-called {\em inverse triangle inequality} at this
point defined by
\begin{equation}\label{A5'}
\sum_{i\in I(x,a(t),b(t))}\lm_i\n a_i(t)\en\le\gg\n\sum_{i\in
I(x,a(t),b(t))}\lm_ia_i(t)\en\;\mbox{ for all }\lm_i\ge 0
\end{equation}
with some constant $\gg>0$; see \cite{ve} for more discussions.

Now we are ready to present the aforementioned {\em well-posedness} result for
the sweeping system \eqref{Problem}, \eqref{C}.\vspace*{-0.1in}

\begin{theorem}[\bf well-posedness of the controlled sweeping
process]\label{Th1} Let all the assumptions in {\rm(H1)}--{\rm(H4)} be
satisfied, and let $(a(\cdot),b(\cdot))\in W^{1,2}([0,T];\R^{mn}\times\R^m)$ and
$u(\cdot)\in L^2([0,T];\R^d)$ be fixed control actions in \eqref{Problem} and
\eqref{C}. Then the sweeping differential inclusion \eqref{Problem} admits the
unique solution $x(\cdot)\in W^{1,2}([0,T];\R^n)$ generated by the control
triple $(a(\cdot),b(\cdot),u(\cdot))$.
\end{theorem}\vspace*{-0.15in}
{\bf Proof.} Following \cite{Tol1}, it is said that a set-valued mapping
$C\colon[0,T]\to\R^n$ is {\em $r$-uniformly lower semicontinuous from the right}
if there exists a family $\mathcal V:=\nn v_r|\;r\ge 0\hnn\subset
W^{1,2}([0,T];\R^n)$ such that for any $r\ge0$, any $s,t\in[0,T]$ with $s\le t$,
and any $x\in\R^n$ with and $\|x\|\le r$ we have the inequality
\begin{equation*}
d\big(x;C(t)\big)\le d\big(x;C(s)\big)+|v_r(t)-v_r(s)|.
\end{equation*}
Let us show that assumption (H4) implies that the polyhedral mapping $C(\cdot)$
defined in \eqref{C} is $r$-uniformly lower semicontinuous from the right. To
proceed, define the function $\phi\colon[0,T]\times\R^n\to\R$ by
\begin{equation}\label{phi}
\phi(t,x):=\disp\max_{1\le i\le m}\big\{\la a_i(t),x\ra-b_i(t)\big\},\quad
t\in[0,T],\;x\in\R^n,
\end{equation}
which gives us the representation $C(t)=\nn x\in\R^n|\;\phi(t,x)\le 0\hnn$ of
the set $C(t)$ from \eqref{C} for each $x\in\R^n$. Let us show that function
\eqref{phi} satisfies the hypothesis $H(\phi)$ formulated in \cite[p.\
297]{Tol1}. Indeed, the convexity of $x\mapsto\phi(t,x)$ and estimates in
\cite[(4.2)]{Tol1} imposed in $H(\phi)$ follow directly from the construction of
$\phi$. Furthermore, we deduce from \eqref{Slater} that the required condition
\cite[(2)]{Tol1} in $H(\phi)$ is also satisfied. To verify $H(\phi)$, it remains
checking the validity of \cite[(4.1)]{Tol1}. Since $b_i(\cdot)\in
W^{1,2}([0,T];\R)$, we clearly have that $\disp\max_{1\le i\le
m}\disp\max_{t\in[0,T]}|b_i(t)|<\infty$. Moreover, it follows from \eqref{phi}
for all $x\in\R^n$ and all $t,s\in[0,T]$ that
\begin{equation*}
|\phi(t,x)-\phi(s,x)|\le\max_{1\le i\le m}\n a_i(t)-a_i(s)\en\cdot\n
x\en+\max_{1\le i\le m}|b_i(t)-b_i(s)|.
\end{equation*}
Taking $v_i(\cdot)$ from assumption (H4), denote further
\begin{equation*}
\xi_r(t):=\disp\int^t_0\(\max_{1\le i\le m}\left|\dot v_i(\t)\right|+\max_{1\le
i\le m}|\dot b_i(\t)|\)\;d\t,\quad r\ge 0.
\end{equation*}
Then $\xi_r(\cdot)\subset W^{1,2}([0,T];\R)$ for all $r\ge 0$, and we have from
\eqref{v-i} that
\begin{equation*}
|\phi(t,x)-\phi(s,x)|\le|\xi_r(t)-\xi_r(s)|\;\mbox{ whenever }\;\|x\|\le
r,\;t,s\in[0,T],
\end{equation*}
which completes the verification of all the assumptions in $H(\phi)$ of
\cite{Tol1}. Employing now \cite[Theorem~4.1]{Tol1} verifies that our polyhedral
mapping $C(\cdot)$ is $r$-uniformly lower semicontinuous from the right on
$[0,T]$. Finally, the existence and uniqueness result claimed in the theorem
follow from \cite[Lemma~3.1 and Theorem~4.1]{Tol1}. $\h$\vspace*{-0.2in}

\section{Existence of Optimal Solutions and Relaxation}\label{exist}
\setcounter{equation}{0}\vspace*{-0.1in}

This section addresses the existence issue for (global) {\em optimal solutions}
to the sweeping control problem $(P)$. Then we define an appropriate notion of
{\em local minimizers} to $(P)$ and discuss its {\em relaxed} counterpart.

Before establishing the existence of optimal solutions to $(P)$ in the
aforementioned class of feasible solutions, let us reformulate the sweeping
differential inclusion \eqref{Problem} in a more convenient way. Consider the
image of the control set $U$ under the perturbation mapping
$g\colon\R^n\times\R^d\to\R^m$ defined by
\begin{equation*}
g(x,U):=\big\{v\in\R^m\big|\;v=g(x,u)\;\mbox{ for some }\;u\in
U\big\},\;x\in\R^n.
\end{equation*}
Then the sweeping inclusion \eqref{Problem} with the moving set \eqref{C} can be
equivalently represented as
\begin{equation}\label{e:3.2}
-\dot{x}(t)\in N\big(x(t);C(t)\big)-g\big(x(t),U\big)\;\mbox{
a.e.}\;t\in[0,T],\;x_0\in C(0).
\end{equation}
More rigorously, this equivalence takes into account standard {\em measurable
selection} results ensuring that for any measurable velocity function satisfying
$v(t)\in g(x(t);U)$ for a.e.\ $t\in[0,T]$ there exists a measurable control
$u(t)\in U$ such that $v(t)=g(x(t),u(t))$ a.e.\ on $[0,T]$. This is surely the
case in our setting; see, e.g., \cite[Chapter~14]{rw} for more details and
references.

Now we are ready to obtain the existence theorem for optimal solutions to $(P)$
under certain additional {\em convexity} assumptions with respect to velocities.
For simplicity we suppose here that the integrand $\ell$ does not depend on the
control variable $u$.  If it does, we have to impose the convexity of an
extended velocity set that includes the integrand component.\vspace*{-0.1in}

\begin{theorem}[\bf existence of optimal solutions to controlled sweeping
processes]\label{Th3} Let $(P)$ be the optimal control problem formulated in
Section~{\rm\ref{intro}} with the equivalent form \eqref{e:3.2} of the sweeping
differential inclusion over all the $W^{1,2}([0,T];\R^n)\times
W^{1,2}([0,T];\R^{mn})\times W^{1,2}([0,T];\R^m)\times L^2([0,T];\R^d)$
quadruples $(x(\cdot),a(\cdot),b(\cdot),u(\cdot))$. In addition to the standing
assumptions {\rm(H1)--(H5)}, suppose that the integrand $\ell$ in
\eqref{minimize} does not depend on the $u$-variable while being convex with
respect to the velocity variables $(\dot x,\dot a,\dot b)$. Suppose furthermore
that along a minimizing sequence of $\(x^k(\cdot),a^k(\cdot),b^k(\cdot),u^k(\cdot)\)$
as $k\in\N$ we have that
$\ell(t,\cdot)$ is majorized by a summable function, that
$\{(x^k(\cdot),a^k(\cdot),b^k(\cdot))\}$ is bounded in
$W^{1,2}([0,T];\R^n\times\R^{mn}\times\R^m)$, and that the set $g(x^k(t);U)$ is
convex for all $t\in[0,T]$. Then $(P)$ admits an
optimal solution in $W^{1,2}([0,T];\R^{n+mn+m})\times L^2([0,T];\R^d)$.
\end{theorem}\vspace*{-0.1in}
{\bf Proof.} Since the set of feasible solutions to problem $(P)$ is nonempty by
Theorem~\ref{Th1}, we can take the minimizing sequence of quadruples
$(x^k(\cdot),a^k(\cdot),b^k(\cdot),u^k(\cdot))$ in $(P)$ from the formulation of
the theorem. It follows from the boundedness of
$\{x^k(\cdot),(a^k(\cdot),b^k(\cdot))\}$ in
$W^{1,2}([0,T];\R^n\times\R^{mn}\times\R^m)$ and the weak compactness of the
dual ball in $L^2([0,T];\R^n\times\R^{mn}\times\R^m)$ that $\dot x^k(\cdot)\to
v^x(\cdot)$, $\dot a^k(\cdot)\to v^a(\cdot)$, and $\dot b^k(\cdot)\to
v^b(\cdot)$ weakly in $L^2([0,T];\R^n)$, $L^2([0,T];\R^{mn})$, and
$L^2([0,T];\R^m)$ along subsequences (without relabeling) for some functions
$v^x(\cdot)$, $v^a(\cdot)$, and $v^b(\cdot)$ from the corresponding spaces.
Employing Mazur's weak closure theorem, we conclude that there are sequences of
convex combinations of $\dot x^k(\cdot)$, $\dot a^k(\cdot)$, and $\dot
b^k(\cdot)$, which strongly converge in the corresponding spaces to
$v^x(\cdot)$, $v^a(\cdot)$, and
$v^b(\cdot)$, respectively. Furthermore, standard real analysis tells us that
there exists a subsequence of these convex combinations (no relabeling again),
which converges to $(v^x(\cdot),v^a(\cdot),v^b(\cdot)$ as $k\to\infty$ a.e.\
pointwise on $[0,T]$. Define now $\ox(\cdot)\in W^{1,2}([0,T];\R^n)$,
$\oa(\cdot)\in W^{1,2}([0,T];\R^{mn})$, and $\ob(\cdot)\in W^{1,2}([0,T];\R^m)$
by
\begin{equation*}
\ox(t):=x_0+\int^t_0v^a(s)ds,\;\oa(t):=a_0+\int^t_0v^a(s)ds,\;\;\mbox{and}
\;\;\ob(t):=b_0+\int^t_0v^b(s)ds,\quad t\in[0,T],
\end{equation*}
and observe that they satisfy the pointwise constraints in \eqref{h1.5} and
\eqref{e:8}. Furthermore, it follows from the closedness and convexity of the
normal cone \eqref{NC} to the moving convex polyhedral set $C(t)$ in \eqref{C}
and the assumed convexity of the compact sets $g(x^k(t),U)$ on $[0,T]$ that the
right-hand site velocity set in \eqref{F-rep} is convex along the selected
minimizing sequence, and we have
\begin{equation*}
\dot{\ox}(t)\in-N\big(\ox(t);\Bar C(t)\big)+g\big(\ox(t),U\big)\;\textrm{ a.e.
}\;t\in[0,T],\quad\ox(0)=x_0\in\Bar C(0)
\end{equation*}
for the limiting trajectory $\ox(\cdot)$ with $\ox(t)\in\Bar
C(t):=\{x\in\R^n|\;\la\oa_i(t),x\ra\le\ob_i(t),\;i=1,\ldots,m\}$ on $[0,T]$.
Employing now the aforementioned measurable selection allows us to find a
measurable control $\ou(\cdot)$ such $\ou(t)\in U$ and
\begin{equation*}
\dot{\ox}(t)\in-N\big(\ox(t);\Bar C(t)\big)+g\big(\ox(t),\ou(t)\big)\;\textrm{
a.e. }\;t\in[0,T].
\end{equation*}
It remains to show that the limiting quadruple
$(\ox(\cdot),\oa(\cdot),\ob(\cdot),\ou(\cdot))$, which is proved to be feasible
for $(P)$, is an optimal solution to this problem. This is a consequence of the
inequality
\begin{equation}\label{cost-lim}
J[\ox,\oa,\ob,\ou]\le\liminf_{k\to\infty}J[x^k,a^k,b^k,u^k]
\end{equation}
for the cost functional \eqref{minimize}. To verify \eqref{cost-lim}, we use the
assumptions in (H5) ensuring the application of the Lebesgue dominated
convergence theorem together with the imposed convexity of integrand with
respect to $(\dot x,\dot a,\dot b)$. This allows us to apply the classical lower
semicontinuity result for integral functionals with respect to the weak topology
in $L^2$. Observe that there is no need to care about the convergence with
respect to $u$-controls in our setting due to the independence of the integral
$\ell$ on the $u$-component. Thus the proof is complete. $\h$

Justifying the existence of global optimal solutions to the controlled sweeping
process under $(P)$, recall that our goal is the derivation of necessary
optimality conditions for suitable {\em local} minimizers of $(P)$ by employing
the method of discrete approximations. An appropriate concept from this
viewpoint goes back to {\em intermediate} local minimizers introduced in
\cite{m95} for Lipschitzian differential inclusions that occupies an
intermediate position between the conventional notions of weak and strong
minimizers in dynamic optimization while covering the latter; see the books
\cite{m-book,v} and the references therein for more details on this notion for
Lipschitzian inclusions. In the case of our problem $(P)$, a natural
implementation of this concept reads as follows.\vspace*{-0.1in}

\begin{definition}[intermediate local minimizers for sweeping optimal
control]\label{Def1}
Let $(\ox(\cdot),\oa(\cdot),\ob(\cdot),\ou(\cdot))$ be a feasible solution to
problem $(P)$ under the standing assumptions made. We say that
$(\ox(\cdot),\oa(\cdot),\ob(\cdot),\ou(\cdot))$ is an {\sc intermediate local
minimizer} {\rm(i.l.m.)} for $(P)$ if
$(\ox(\cdot),\oa(\cdot),\ob(\cdot),\ou(\cdot))\in W^{1,2}([0,T];\R^n)\times
W^{1,2}([0,T];\R^{mn})\times W^{1,2}([0,T];\R^m)\times L^2([0,T];\R^d)$ and
there exists $\ve>0$ such that
\begin{equation*}
J[\ox,\oa,\ob,\ou]\le J[x,a,b,u]
\end{equation*}
for any feasible solutions $(x(\cdot),a(\cdot),b(\cdot),u(\cdot))$ to $(P)$
satisfying
\begin{equation}\label{ilm}
\big\|x(\cdot)-\ox(\cdot)\big\|_{W^{1,2}}+\big\|\big(a(\cdot),
b(\cdot)\big)-\big(\oa(\cdot),\ob(\cdot)\big)\big\|_{W^{1,2}}+\|
u(\cdot)-\ou(\cdot)\|_{L^2}\le\ve.
\end{equation}
\end{definition}\vspace*{-0.1in}
If the term $\|x(\cdot)-\ox(\cdot)\|_{W^{1,2}}$ in \eqref{ilm} is replaced by
$\|x(\cdot)-\ox(\cdot)\|_{\cal C}$, the norm in the space of continuous
functions ${\cal C}([0,T];\R^n)$, we speak about {\em strong local minimizers}
for $(P)$. It is clear that any strong local minimizer for $(P)$ is an
intermediate one, but not vice versa as can be confirmed by
examples.\vspace*{-0.03in}

To implement our approach to study local minimizers of $(P)$, we need a certain
{\em relaxation stability} of the i.l.m.\ under consideration. The idea of
relaxation of variational problems, related to convexification with respect to
derivative variables, goes back to Bogolyubov and Young for the classical
calculus of variations and to Gamkrelidze and Warga for optimal control problems
governed by ordinary differential equations; see, e.g., the books
\cite{m-book,v} for more discussions and references, where relaxation of control
problems for Lipschitzian differential inclusions were also investigated and
discussed in detail. Relaxation results for non-Lipschitzian differential
inclusions were more recently developed in \cite{dfm,et,Tol}.\vspace*{-0.05in}

To proceed in the case of our optimal control problem $(P)$, consider vectors
$x:=(x_1,\ldots,x_n)\in\R^n$, $a:=(a_1,\ldots,a_m)\in\R^{mn}$,
$b:=(b_1,\ldots,b_m)\in\R^m$, and $u:=(u_1,\ldots,u_d)\in\R^d$, and then define
the set-valued mapping $F\colon\R^n\times\R^{mn}\times\R^m\times\R^d\tto\R^n$ by
\begin{equation}\label{F}
F(x,a,b,u):=N\big(x;C(a,b)\big)-g(x,u),
\end{equation}
where $N(x;C(a,b))$ is taken from in \eqref{NC}, and where $C(a,b):=\nn
x\in\R^n\;|\;\la a_i,x\ra\le b_i,\;i=1,\ldots,m\hnn$. It is not hard to see that
$F$ admits the following explicit representation:
\begin{equation}\label{F-rep}
F(x,a,b,u)=\Big\{\sum_{i\in I(x,a,b)}\eta_ia_i\Big|\;\eta_i\ge 0\Big\}-g(x,u)
\end{equation}
via the active index set \eqref{e:AI} at $x\in C(a,b)$. Let $\ell_F(t,a,b,u,\dot
x,\dot a,\dot b)$ be the restriction of the integral $\ell$ on the set
$F(x,a,b,u)$ with $\ell_F(t,a,b,u,\dot x,\dot a,\dot b):=\emp$ if $\dot x\notin
F(x,a,b,u)$. Denoting by $\Hat\ell_F$ the {\em convexification} of the integrand
(i.e., the largest l.s.c.\ convex function majorized by
$\ell(t,x,a,b,\cdot,\cdot,\cdot,\cdot)$) with respect to the velocity variables
$(\dot x,\dot a,\dot b)$ as well as to the control one $u$ on the convex hull
$\co U$, define the {\em relaxed optimal control problem} $(R)$ by:
\begin{equation}\label{R}
\textrm{minimize}\hspace{0.1in}\hat
J[x,a,b,u]:=\vph\big(x(T)\big)+\int_0^T\hat\ell_F\big(t,x(t),a(t),b(t),u(t),\dot
{x}(t),\dot{a}(t),\dot{b}(t)\big)dt
\end{equation}
over quadruples $(x(\cdot),a(\cdot),b(\cdot),u(\cdot))\in
W^{1,2}([0,T];\R^n)\times W^{1,2}([0,T];\R^{mn})\times W^{1,2}([0,T];\R^m)\times
L^2([0,T];\R^d)$ satisfying \eqref{h1.5} and giving a finite value of the
extended running cost in \eqref{R}. All such quadruples are said to be {\em
feasible} to $(R)$. It follows from \eqref{R} and the construction of
$\hat\ell_F$ with $F$ taken from \eqref{F} that $u(t)\in\co U$ for a.e.\
$t\in[0,T]$, and that $x(\cdot)$ is a trajectory of the convexified differential
inclusion
\begin{equation}\label{conv}
-\dot{x}(t)\in N\big(x(t);C(t)\big)-\co g\big(x(t),U\big)\;\mbox{
a.e.}\;t\in[0,T],\;x_0\in C(0)
\end{equation}
with $\la a_i(t),x(t)\ra\le b_i(t)$ for $i=1,\ldots,m$ and all $t\in[0,T]$. Now
we introduce a new notion of {\em relaxed} intermediate local minimizers for
$(P)$; cf.\ \cite{m95} for Lipschitzian differential inclusions and \cite{cmn}
for a version of problem $(P)$ with $\ell\equiv 0$ and an uncontrolled
polyhedron $C(t)\equiv C$.\vspace*{-0.1in}

\begin{definition}[\bf relaxed intermediate local minimizers]\label{Def2} We say
that $(\ox(\cdot),\oa(\cdot),\ob(\cdot),\ou(\cdot))$ is a {\sc relaxed
intermediate local minimizer} $($r.i.l.m.$)$ for problem $(P)$ if it is feasible
for $(P)$ and there exists $\ve>0$ such that
\begin{equation*}
J[\ox,\oa,\ob,\ou]=\hat J[\ox,\oa,\ob,\ou]\le\hat J[x,a,b,u]
\end{equation*}
whenever a feasible quadruple $(x(\cdot),a(\cdot),b(\cdot),u(\cdot))$ for $(R)$
satisfies \eqref{ilm}.
\end{definition}\vspace*{-0.05in}

It follows from Definitions~\ref{Def1} and \ref{Def2} in view of the
constructions in \ref{R} and \eqref{conv} that any i.l.m.\ of $(P)$ is also its
r.i.l.m.\ provided that the sets $U$ and $g(x(t);U)$ are convex and the
integrand $\ell(t,x(t),a(t),b(t),\cdot,\cdot,\cdot,\cdot)$ is convex along
feasible solutions to $(P)$. The well-recognized beauty of relaxation procedures
in variational and control problems is that they keep global or local optimal
values of cost functionals under relaxation in important situations {\em without
any convexity assumptions}. It is strongly related to deep measure-theoretical
results of the Lyapunov-Aumann type ensuring the automatic convexity of
integrals of arbitrary set-valued mappings over nonatomic measures. In
particular, it has been realized in this way that every strong local minimizer
in control problems for Lipschitzian differential inclusions with no constraint
of right ends of trajectories is always a relaxed one; see, e.g.,
\cite{m-book,v}. Similar results for
controlled sweeping processes of different types were obtained in
\cite[Theorem~2]{et} and \cite[Theorem~4.2]{Tol}. We {\em conjecture} that
modifying the proofs of the aforementioned theorems lead us to the fact that any
{\em strong} local minimizer of the nonconvex sweeping control problem $(P)$ is
a relaxed strong local minimizer of this problem under the imposed standing
assumptions in (H1)--(H5) with the replacement of the lower semicontinuity of
$\ph$ and $\ell$ in (H5) by their continuity.\vspace*{-0.2in}

\section{Strong Discrete Approximation of Feasible Solutions}\label{disc}
\setcounter{equation}{0}\vspace*{-0.1in}

In this section we start our detailed development of the method of discrete
approximations to study the sweeping optimal control problem $(P)$ formulated in
Section~\ref{intro}. In fact, this section does not concern the optimization
part of $(P)$ while dealing only with constructive approximations of {\em
feasible solutions}. Our main goal here is to show that the standing assumptions
imposed allow us to {\em strongly} approximate {\em any} feasible solution to
$(P)$ by feasible solutions to discrete-time problems extended to the
continuous-time interval. The result established below significantly improves
similar ones obtained in \cite{cm1,cm3,chhm3} for particular types of sweeping
control problems, and so its proof is more involved in comparison with those
given in \cite{cm1,cm3,chhm3}. Note that another discrete approximation scheme
was developed in \cite{cmn} for problem $(P)$ with $\ell\equiv 0$ and an
uncontrolled polyhedral convex set $C(t)\equiv C$.

To proceed, for each $k\in\N$ define the discrete partition of $[0,T]$ by
\begin{equation}\label{e:DP}
\Delta_k:=\left\{0=t^k_0<t^k_1<\ldots<t^k_{\nu(k)-1}<t^k_{\nu(k)}=T\right\}
\;\;\mbox{with}\;\;h^k_j:=t^k_{j+1}-t^k_j\le
\dfrac{\Tilde\nu}{\nu(k)}\;\;\mbox{for}\;\;j=0,\ldots,\nu(k)-1,
\end{equation}
where $\nu=\nu(k)\ge k$, and where $\Tilde\nu>0$ is some constant.

Here is a major approximation result, which certainly is of its own interest
(also from a numerical viewpoint), while being important for the subsequent
developments of this paper and its continuation in \cite{ccmn}.\vspace*{-0.05in}

\begin{theorem}[\bf strong discrete approximation of feasible sweeping
solutions]\label{Th2} Under the standing assumptions in {\rm(H1)--(H4)}, fix
any feasible solution $(\ox(\cdot),\oa(\cdot),\ob(\cdot),\ou(\cdot))$ to $(P)$
such that the functions $\dot\ox(\cdot),\dot\oa(\cdot),\dot\ob(\cdot)$ and
$\ou(\cdot)$ are of bounded variation on $[0,T]$, i.e.,
\begin{equation}\label{e:var}
\max\nn{\rm var}\big(\dot\ox(\cdot);[0,T]\big),{\rm
var}\big(\dot\oa(\cdot);[0,T]\big),{\rm var}\big(\dot\ob(\cdot);[0,T]\big),{\rm
var}\big(\ou(\cdot);[0,T]\big)\hnn\le K
\end{equation}
for some constant $K>0$. Then there exist partitions $\Delta_k$, $k=1,\ldots$,
as in \eqref{e:DP} together with sequences of piecewise linear
functions $(x^k(t),a^k(t),b^k(t))$ and piecewise constant functions $u^k(\cdot)$
on $[0,T]$, as well as a sequence of positive numbers $\dd_k$ converging to zero
such that $(x^k(0),a^k(0),b^k(0))=(x_0,a_0,b_0)$ for all $k\in\N$, and we have
the relationships:
\begin{equation}\label{e:a-dc}
1-\dd_k\le\left\|a^k_i(t^k_j)\right\|\le 1+\dd_k\;\mbox{ for all
}\;t^k_j\in\Delta_k,\quad i=1,\ldots,m,
\end{equation}
\begin{equation*}
x^k(t)=x^k(t^k_j)+(t-t^k_j)v^k_j,\;\;t^k_j\le t\le
t^k_{j+1}\;\;\mbox{with}\;\;-v^k_j\in
F\big(x^k(t^k_j),a^k(t^k_j),b^k(t^k_j),u^k(t^k_j)\big)
\end{equation*}
for $j=0,\ldots,\nu(k)-1$ together with the convergence
$\nn(x^k(\cdot),a^k(\cdot),b^k(\cdot))\hnn\to(\ox(\cdot),\oa(\cdot),\ob(\cdot))$
in the $W^{1,2}$-norm topology on $[0,T]$, and $\nn u^k(\cdot)\hnn\to\ou(\cdot)$
in the $L^2$-norm topology on $[0,T]$ as $k\to\infty$.
\end{theorem}\vspace*{-0.1in}
{\bf Proof.} We split the proof into the following four major steps.\\[1ex]
First of all, we remark that, thanks to our assumption \eqref{e:var}, we can suppose $\dot{\bar{x}}, \dot{\bar{a}},\dot{\bar{b}}$ and
$\bar{u}$ to be defined everywhere.\\[1ex]
{\bf Step~1:} {\em Constructing $\(u^k(\cdot),a^k(\cdot)\)$ to approximate
$\(\ou(\cdot),\oa(\cdot)\)$}. Since step functions are dense in $L^2[0,T]$,
there are sequences of step functions $\nn
u^k(\cdot)\hnn=\nn\(u^k_1(\cdot),\ldots,u^k_d(\cdot)\)\hnn$ and
$\nn\al^k(\cdot)\hnn=\nn\(\al^k_1(\cdot),\ldots,\al^k_m(\cdot)\)\hnn$ with
\begin{equation}\label{mu-k}
\mu_k:=\max\nn\int^T_0\n u^k(t)-\ou(t)\en^2
dt,\int^T_0\n\al^k(t)-\dot\oa(t)\en^2 dt\hnn\to 0\;\;\mbox{as}\;\;k\to\infty.
\end{equation}
Furthermore, for each $k\in\N$ we find a partition $\Delta_k$ of the interval
$[0,T]$ from \eqref{e:DP} for which the step functions $\nn u^k(\cdot)\hnn$ and
$\nn\al^k(\cdot)\hnn$ are constant on the subintervals $[t_j,t_{j+1})$ for
$j=0,\ldots,\nu(k)-1$.
This gives us the strong convergence of
$\nn\(u^k(\cdot),\al^k(\cdot)\)\hnn$ to $\(\ou(\cdot),\dot{\bar{x}}(\cdot)\)$ in
$L^2([0,T])$ as $k\to\infty$. Since the intervals $(t_j,t_{j+1})$ are not prescribed \textit{a priori} and $\bar{x}$
is a Caratheodory solution of \eqref{e:3.2}, up to possibly increasing the number of intervals of the partition we
can suppose without loss of generality that the differential inclusion \eqref{e:3.2} is satisfied at all endpoints of
$\delta_k$ that are contained in the open interval $(0,T)$.
Next we define the functions $a^k(\cdot)$ by
\begin{equation}\label{e:ak}
a^k(t):=a_0+\int^t_0\al^k(s)ds,\quad t\in[0,T].
\end{equation}
It tells us that each $a^k(\cdot)$ is piecewise linear on $[0,T]$, since its
derivative $\dot a^k(\cdot)=\al^k(\cdot)$ is piecewise constant on $[0,T]$. By
\eqref{mu-k} we have the strong convergence in $L^2([0,T])$ of
$\nn\dot a^k(\cdot)\hnn$ to $\dot\oa(\cdot)$. Moreover, it follows from
\eqref{e:ak} and the classical H\"{o}lder inequality that
\begin{equation}\label{e:a-est}
\left
|a^k_{ip}(t)-\oa_{ip}(t)\right|^2=\left|\int^t_0\[\al^k_{ip}(s)-\oal_{ip}(s)\]
ds\right|^2 \le\[\int^T_0\left|\al^k_{ip}(s)-\oal_{ip}(s)\right|^2
ds\]T\le\mu_kT
\end{equation}
for all $t\in[0,T]$, $i=1,\ldots,m$, and  each component index $p=1,\ldots,n$.
Hence the sequence of functions $a^k(\cdot)$ converges strongly to $\oa(\cdot)$
in $W^{1,2}([0,T])$ and satisfies the estimates in \eqref{e:a-dc} with
\begin{equation}\label{e:def-deltak}
\dd_k:=\sqrt{n\mu_kT}.
\end{equation}
{\bf Step~2:} {\it Constructing $\(x^k(\cdot),b^k(\cdot)\)$ to approximate
$\(\ox(\cdot),\ob(\cdot)\)$.} While proceeding recurrently, fix any
$j\in\{0,\ldots,\nu(k-1)\}$, suppose that the pairs $(x^k_j,b^k_j)$ are known
for all $j=0,\ldots,\nu(k)-1$, and then construct the pair
$(x^k_{j+1},b^k_{j+1})$. Define the numbers
\begin{equation}\label{e:bjk}
b^k_{ij}:=\la
a^k_i(t_j),x^k_j\ra+\ob_i(t_j)-\la\oa_i(t_j),\ox(t_j)\ra\;\;\mbox{for
all}\;\;i=1,\ldots,m,
\end{equation}
\begin{equation}\label{e:bk}
b^k_i(0):=b_{i0},\;\;\mbox{and}\;\;b^k_i(t):=b^k_{ij}+\dfrac{t-t_j}{h^k_j}
\(b^k_{i,j+1}-b^k_{ij}\)\;\;\mbox{for all}\;\;
t\in[t_j,t_{j+1}]\;\;\mbox{and}\;\;i=1,\ldots,m.
\end{equation}
It gives us $b^k_{ij}-\la
a_i^k(t_j),x^k_j\ra=\ob_i(t_j)-\la\oa_i(t_j),\ox(t_j)\ra$, and hence
\begin{equation}\label{index}
I(x^k_j,a^k(t_j),b^k_j)=I\big(\ox(t_j),\oa(t_j),\ob(t_j)\big)\;\;\mbox{for
all}\;\;j=0,\ldots,\nu(k).
\end{equation}
It follows from the validity of $-\dot{\ox}(t)=F(\ox(t),\oa(t),\ob(t),\ou(t))$
for a.e.\ $t\in[0,T]$ including the mesh points of $\Delta_k$ with $F$ given in
\eqref{F}, the measurability of the set-valued mapping $t\mapsto
F(\ox(t),\oa(t),\ob(t),\ou(t))$ on $[0,T]$ due to \cite[Theorem~14.26]{rw} with
the representation of $F$ in \eqref{F-rep}, and the measurable selection result
from \cite[Corollary~14.6]{rw} that there exist nonnegative measurable functions
$\eta_i(\cdot)$ on $[0,T]$ as $i=1,\ldots,m$ ensuring the equality
\begin{equation*}
-\dot\ox(t)=\sum_{i\in
I(\ox(t),\oa(t),\ob(t))}\eta_i(t)\oa_i(t)-g(\ox(t),\ou(t)) \;\;\mbox{for
a.e.}\;\;t\in[0,T].
\end{equation*}
Define now the vectors $v^k_j$ for all indices $j=0,\ldots,\nu(k)$ by
\begin{equation}\label{e:dv}
-v^k_j:=\sum_{i\in
I(\ox(t_j),\oa(t_j),\ob(t_j))}\eta_i(t_j)a^k_i(t_j)-g\big(x^k_j,
u^k(t_j)\big)=\sum_{i\in
I(x^k_j,a^k(t_j),b^k_j)}\eta_i(t_j)a^k_i(t_j)-g\big(x^k_j,u^k(t_j)\big),
\end{equation}
where the second equality comes from \eqref{index}. It is obvious that
$-v^k_j\in F(x^k_j,a^k(t_j),b^k(t_j),u^k(t_j))$ for such indices $j$.
Since $\dot x(\cdot)$ is of bounded variation on $[0,T]$, we have by
\eqref{e:var} that
\begin{equation*}
\n\dot\ox(t)-\dot\ox(0)\en\le\n\dot\ox(t)-\dot\ox(0)\en+\n\dot\ox(T)-\dot\ox(t)
\en\le{\rm var}\big(\dot\ox(\cdot);[0,T]\big)\le K.
\end{equation*}
which in turn yields the estimate
\begin{equation*}
\n\dot\ox(t)\en\le\n\dot\ox(0)\en+K:=M^x_1
\end{equation*}
for a.e.\ $t\in[0,T]$ including the mesh points of $\Delta_k$. Using the inverse
triangle inequality \eqref{A5'} implies that
$$
\begin{aligned}
\eta_i(t)=\eta_i(t)\n\oa_i(t)\en&\le\sum_{i\in
I(\ox(t),\oa(t),\ob(t))}\eta_i(t)\n\oa_i(t)\en\le\gg\n\sum_{i\in
I(\ox(t),\oa(t),\ob(t))}\eta_i(t)\oa_i(t)\en\\
&\le\gg\n\dot\ox(t)\en+\gg\n g\big(\ox(t),\ou(t)\big)\en\le\gg M^x_1+\gg
M\(1+\|\ox(t)\|\)\\
&\le\gg M^x_1+\gg M\(1+\max_{t\in[0,T]}\|\ox(t)\|\)=:M^x_2
\end{aligned}
$$
for a.e.\ $t\in[0,T]$ and for all $i\in I(\ox(t),\oa(t),\ob(t))$. By
\eqref{e:gL} it yields the estimates
\begin{equation}\label{e:v-est}
\begin{aligned}
\n v^k_j-\dot{\ox}(t_j)\en&\le\sum_{i\in
I(\ox(t_j),\oa(t_j),\ob(t_j))}\eta_i(t_j)\n\oa_i(t_j)-a^k_i(t_j)\en+\n
g\big(\ox(t_j),\ou(t_j)\big)-g\big(x^k_j,u^k(t_j)\big)\en\\
&\le
M^x_2\sum_{i=1}
^m\n\oa_i(t_j)-a^k_i(t_j)\en+L\(\n\ox(t_j)-x^k_j\en+\n\ou(t_j)-u^k(t_j)\en\).
\end{aligned}
\end{equation}
Letting now $x^k_{j+1}:=x^k_j+h^k_jv^k_j$, we define the arcs $x^k(t)$ on
$[0,T]$
\begin{equation}
\label{e:xk}
x^k(t):=x^k_j+\dfrac{t-t_j}{h^k_j}\(x^k_{j+1}-x^k_j\)=x^k_j+(t-t_j)v^k_j,
\;\;\mbox{for}\;\;t\in[t_j,t_{j+1}]
\end{equation}
and thus complete the construction of the pairs $\(x^k(\cdot),b^k(\cdot)\)$ in
this step.\\[1ex]
{\bf Step~3:} {\em Verifying the strong $W^{1,2}$-convergence of $x^k(\cdot)$ to
$\ox(\cdot)$ on $[0,T]$}. For each index $j=0,\ldots,\nu(k)-1$ and
$i=1,\ldots,m$ denote the functions on $[t_j,t_{j+1})$ by
$$
f^x_j(s):=\n\dot\ox(t_j)-\dot\ox(s)\en,\;f^u_j(s):=\n\ou(t_j)-\ou(s)\en,\;f^a_{
ij}(s):=\n\oal_i(t_j)-\oal_i(s)\en,\;f^b_{ij}(s):=\n\obe_i(t_j)-\obe_i(s)\en
$$
and then select $s^x_j,s^u_j,s^a_{ij},s^b_{ij}$ from the subintervals
$[t_j,t_{j+1})$ such that
\begin{equation}\label{BV}
\begin{cases}
\disp\sup_{s\in[t_j,t_{j+1}]}f^x_j(s)\le\n\dot\ox(t_j)-\dot\ox(s^x_j)\en+2^{-k},
\\
\disp\sup_{s\in[t_j,t_{j+1}]}f^u_j(s)\le\n\ou(t_j)-\ou(s^u_j)\en+2^{-k},\\
\disp\sup_{s\in[t_j,t_{j+1}]}f^a_{ij}(s)\le\n\oal(t_j)-\oal(s^a_{ij})\en+2^{-k},
\\
\disp\sup_{s\in[t_j,t_{j+1}]}f^b_{ij}(s)\le\left|\obe(t_j)-\obe(s^b_{ij}
)\right|+2^{-k}.\\
\end{cases}
\end{equation}
With $h_k:=\disp\max_{0\le j\le\nu(k)-1}\{h^k_j\}$, we get from the above the
following relationships:
\begin{equation}\label{e:x-est}
\begin{aligned}
\n x^k_{j+1}-\ox(t_{j+1})\en &=\n
x^k_j+h^k_jv^k_j-\ox(t_j)-\int^{t_{j+1}}_{t_j}\dot\ox(s)ds\en\le\n
x^k_j-\ox(t_j)\en+\int^{t_{j+1}}_{t_j}\n v^k_j-\dot\ox(s)\en ds\\
&\le\n x^k_j-\ox(t_j)\en+\int^{t_{j+1}}_{t_j}\n v^k_j-\dot\ox(t_j)\en
ds+\int^{t_{j+1}}_{t_j}\n\dot\ox(t_j)-\dot\ox(s)\en ds\\
&\le\n x^k_j-\ox(t_j)\en+h^k_jL\(\n x^k_j-\ox(t_j)\en+\n\ou(t_j)-u^k(t_j)\en\)\\
&+h^k_j
M^x_2\sum_{i=1}^m\n\oa_i(t_j)-a^k_i(t_j)\en+\int^{t_{j+1}}_{t_j}f^x_j(s)ds\\
&\le\(1+Lh^k_j\)\n x^k_j-\ox(t_j)\en+L\int^{t_{j+1}}_{t_j}\n\ou(t_j)-\ou(s)\en
ds+L\int^{t_{j+1}}_{t_j}\n\ou(s)-u^k(s)\en ds\\
&+M^x_2mh^k_j\dd_k+\int^{t_{j+1}}_{t_j}f^x_j(s)ds\\
&\le\(1+Lh_k\)\n x^k_j-\ox(t_j)\en+L\int^{t_{j+1}}_{t_j}\n\ou(s)-u^k(s)\en
ds+L\int^{t_{j+1}}_{t_j}f^u_j(s)ds+\\
&+\int^{t_{j+1}}_{t_j}f^x_j(s)ds+M^x_2mh^k_j\dd_k.
\end{aligned}
\end{equation}
Let $A:=1+Lh_k$, and for each $j=0,\ldots,\nu(k)-1$ denote $\gg_j:=\n
x^k_j-\ox(t_j)\en$ and
$$
\begin{aligned}
\lm_j:=L\int^{t_{j+1}}_{t_j}\n\ou(s)-u^k(s)\en ds+L\int^{t_{j+1}}_{t_j}
f^u_j(s)ds+\int^{t_{j+1}}_{t_j} f^x_j(s)ds+M^x_2mh^k_j\dd_k.
\end{aligned}
$$
Then the final estimate in \eqref{e:x-est} reads as
\begin{equation*}
\gg_{j+1}\le A\gg_j+\lm_j\;\;\mbox{for}\;\;j=0,\ldots,\nu(k)-1,
\end{equation*}
which in turn implies the conditions
\begin{equation*}
\gg_j\le A^j\gg_0+A^{j-1}\lm_0+A^{j-2}\lm_1+\ldots
A^0\lm_j=A^{j-1}\lm_0+A^{j-2}\lm_1+\ldots A^0\lm_j.
\end{equation*}
Since $A^j=\(1+Lh_k\)^j\le\(1+Lh_k\)^{\nu(k)}\le e^{L\Tilde\nu}$, we get
$\gg_j\le e^{L\Tilde\nu}\(\lm_0+\lm_1+\ldots+\lm_j\)\le
e^{L\Tilde\nu}\disp\sum^{\nu(k)-1}_{j=0}\lm_j$.\vspace*{-0.1in}

Let us next estimate the quantity
\begin{equation}\label{suml}
\sum^{\nu(k)-1}_{j=0}\lm_j=\sum^{\nu(k)-1}_{j=0}\[L\int^{t_{j+1}}_{t_j}
\n\ou(s)-u^k(s)\en ds+L\int^{t_{j+1}}_{t_j}f^u_j(s)ds+
\int^{t_{j+1}}_{t_j}f^x_j(s)ds+M^x_2mh^k_j\dd_k\].
\end{equation}
To proceed, we deduce from \eqref{mu-k} and \eqref{BV} that
\begin{eqnarray}\label{e:u-est}
\begin{array}{ll}
\disp\sum^{\nu(k)-1}_{j=0}\int^{t_{j+1}}_{t_j}\n\ou(s)-u^k(s)\en ds &\le\sqrt
T\disp\sqrt{\sum^{\nu(k)-1}_{j=0}\int^{t_{j+1}}_{t_j}\n\ou(s)-u^k(s)\en^2 ds}\\
&=\disp\sqrt T\sqrt{\int^T_0\n u^k(t)-\ou(t)\en^2 dt}\le\sqrt{T\mu_k},
\end{array}
\end{eqnarray}
\begin{equation}\label{e:u-est1}
\begin{aligned}
\sum^{\nu(k)-1}_{j=0}\int^{t_{j+1}}_{t_j}f^u_j(s)ds&\le
h_k\sum^{\nu(k)-1}_{j=0}\(\n\ou(t_j)-\ou(s^u_j)\en+\n\ou(s^u_j)-\ou(t_{j+1}
)\en+2^{-k}\)\\
&\le h_k{\rm var}\big(\ou;[0,T]\big)+h_k\nu(k)2^{-k}\le h_k\mu+\Tilde\nu2^{-k}.
\end{aligned}
\end{equation}\\
Using the same arguments leads us to the inequalities
$$
\begin{aligned}
\sum^{\nu(k)-1}_{j=0}\int^{t_{j+1}}_{t_j}f^x_j(s)ds&\le
h_k\sum^{\nu(k)-1}_{j=0}
\(\n\dot\ox(t_j)-\dot\ox(s^x_j)\en+\n\dot\ox(s^x_j)-\dot\ox(t_{j+1})\en+2^{-k}
\)\\
&\le h_k{\rm var}\big(\dot\ox;[0,T]\big)+h_k\nu(k)2^{-k}\le
h_k\mu+\Tilde\nu2^{-k}.
\end{aligned}
$$
On the other hand, we clearly have that
$$
\sum^{\nu(k)-1}_{j=0}M^x_2mh^k_j\dd_k\le M^x_2m\nu(k)h_k\dd_k\le
M^x_2m\Tilde\nu\dd_k.
$$
Combining all the above brings us to the the desired estimate of the quantity
\eqref{suml} and hence of $\|x^k_j-\ox(t_j)\|$:
$$
\begin{aligned}
\sum^{\nu(k)-1}_{j=0}\lm_j&\le L(\sqrt{T\mu_k}+h_k\mu+\Tilde\nu2^{-k})+h_k\mu
+\Tilde\nu2^{-k}+M^x_2m\Tilde\nu\dd_k\\
&\le\(h_k\mu+\Tilde\nu2^{-k}\)(L+1)+L\sqrt{T\mu_k}+M^x_2m\Tilde\nu\dd_k,
\end{aligned}
$$
\begin{equation}
\label{e:x-est1}
\n
x^k_j-\ox(t_j)\en\le\vartheta_k:=e^{L\Tilde\nu}\[\(h_k\mu+\Tilde\nu2^{-k}
\)(L+1)+L\sqrt{T\mu_k}+M^x_2m\Tilde\nu\dd_k\]
\end{equation}
for all $j=0,\ldots,\nu(k)$. Employing this together with \eqref{e:xk},
\eqref{e:x-est}, and \eqref{e:x-est1} gives us
$$
\begin{aligned}
\n x^k(t)-\ox(t)\en&=\n x^k_j+h_kv^k_j-\ox(t_j)-\int^t_{t_j}\dot\ox(s)ds\en\le\n
x^k_j-\ox(t_j)\en+\int^t_{t_j}\n v^k_j-\dot\ox(s)\en ds\\
&\le\n x^k_j-\ox(t_j)\en+\int^{t_{j+1}}_{t_j}\n v^k_j-\dot\ox(s)\en
ds\le\(1+Lh_k\)\n x^k_j-\ox(t_j)\en+\lm_j\\
&\le\(1+Lh_k\)\vartheta_k+\lm_j\;\mbox{ whenever }\;t\in(t_j,t_{j+1}]\;\mbox{
and }\;j=0,\ldots,k-1,
\end{aligned}
$$
which justifies by $\lm_j\to 0$ the uniform convergence of the sequence $\nn
x^k(\cdot)\hnn$ to $\ox(\cdot)$ as $k\to\infty$.

To verify further the $L^2$-strong convergence of $\nn\dot x^k(\cdot)\hnn$ to
$\dot{\ox}(\cdot)$ on $[0,T]$ as $k\to\infty$, observe first that{\color{blue}, owing to \eqref{e:v-est},}
$$
\begin{aligned}
h^k_j\n v^k_j-\dot{\ox}(t_j)\en^2&\le
h^k_j\[M^x_2\sum_{i=1}
^m\n\oa_i(t_j)-a^k_i(t_j)\en+L\(\n\ox(t_j)-x^k_j\en+\n\ou(t_j)-u^k(t_j)\en\)\]
^2\\
&\le
3(M^x_2)^2h^k_jm\sum_{i=1}
^m\n\oa_i(t_j)-a^k_i(t_j)\en^2+3Lh^k_j\n\ox(t_j)-x^k_j\en^2+3Lh^k_j\n\ou(t_j)-u^
k(t_j)\en^2\\
&\le 3(M^x_2)^2m^2\dd_k^2h_k+3Lh_k\delta^2_k+3Lh^k_j\n\ou(t_j)-u^k(t_j)\en^2
\end{aligned}
$$
for $j=0,\ldots,\nu(k)-1$ and then subsequently derive the estimates
$$
\begin{aligned}
\int^T_0&\n \dot
x^k(t)-\dot\ox(t)\en^2dt=\sum^{\nu(k)-1}_{j=0}\int^{t_{j+1}}_{t_j}\n
v^k_j-\dot\ox(t)\en^2dt\\
&\le\sum^{\nu(k)-1}_{j=0}\int^{t_{j+1}}_{t_j}\(\n
v^k_j-\dot\ox(t_j)\en+\n\dot\ox(t_j)-\dot\ox(t)\en\)^2dt\\
&\le 2\sum^{\nu(k)-1}_{j=0}\int^{t_{j+1}}_{t_j}\n
v^k_j-\dot\ox(t_j)\en^2dt+2\sum^{k-1}_{j=0}\int^{t_{j+1}}_{t_j}
\n\dot\ox(t_j)-\dot\ox(t)\en^2dt\\
&\le
2\sum^{\nu(k)-1}_{j=0}\[
3(M^x_2)^2m^2\dd_k^2h_k+3Lh_k\delta^2_k+3Lh^k_j\n\ou(t_j)-u^k(t_j)\en^2\]+
2\sum^{\nu(k)-1}_{j=0}\int^{t_{j+1}}_{t_j}\[f^x_j(t)\]^2dt\\
&\le
6(M^x_2)^2m^2\dd_k^2h_k\nu(k)+6Lh_k\nu(k)\delta^2_k+6L\sum^{\nu(k)-1}_{j=0}
h^k_j\n\ou(t_j)-u^k(t_j)\en^2+2
\sum^{\nu(k)-1}_{j=0}\int^{t_{j+1}}_{t_j}\[f^x_j(t)\]^2dt\\
&\le
6(M^x_2)^2m^2\dd_k^2\Tilde\nu+6L\Tilde\nu\delta^2_k+6L\sum^{\nu(k)-1}_{j=0}
h^k_j\n\ou(t_j)-u^k(t_j)\en^2+2
\sum^{\nu(k)-1}_{j=0}\int^{t_{j+1}}_{t_j}\[f^x_j(t)\]^2dt.
\end{aligned}
$$
Since the control set $U$ is compact, there exists a number $\Bar M>0$ such that
$\max\nn\n\ou(t)\en,\n u^k(t)\en\hnn\le\Bar M$ for all $t\in[0,T]$. On the other
hand, it follows from \eqref{e:u-est} and \eqref{e:u-est1} that
$$
\begin{aligned}
\sum^{\nu(k)-1}_{j=0}h^k_j\n\ou(t_j)-u^k(t_j)\en^2&\le 2\Bar
M\sum^{\nu(k)-1}_{j=0}\int^{t_{j+1}}_{t_j}\n\ou(t_j)-u^k(t_j)\en dt\\
&\le 2\Bar M\sum^{\nu(k)-1}_{j=0}\int^{t_{j+1}}_{t_j}\n\ou(t_j)-\ou(t)\en
dt+2\Bar M\sum^{\nu(k)-1}_{j=0}\int^{t_{j+1}}_{t_j}\n\ou(t)-u^k(t)\en dt\\
&\le 2\Bar M\sum^{\nu(k)-1}_{j=0}\int^{t_{j+1}}_{t_j}f^u_j(t)dt+2\Bar
M\sqrt{T\mu_k}\le 2\Bar M\(h_k\mu+\Tilde\nu2^{-k}+\sqrt{T\mu_k}\).
\end{aligned}
$$
In addition we get from the constructions and notation above that
$$
\begin{aligned}
\sum^{\nu(k)-1}_{j=0}\int^{t_{j+1}}_{t_j}\[f^x_j(t)\]^2dt&\le\sum^{\nu(k)-1}_{
j=0}\int^{t_{j+1}}_{t_j}\[f^x_j(s_j)+2^{-k}\]^2dt=\sum^{\nu(k)-1}_{j=0}h^k_j\[
f^x_j(s_j)+2^{-k}\]^2\\
&\le 2 h_k\sum^{\nu(k)-1}_{j=0}\nn\[f^x_j(s_j)\]^2+4^{-k}\hnn\le
2h_k\[\sum^{\nu(k)-1}_{j=0}f^x_j(s_j)\]^2+2h_k\nu(k)4^{-k}\\
&\le 2
h_k\[\sum^{\nu(k)-1}_{j=0}
\(\n\dot\ox(t_j)-\dot\ox(s^x_j)\en+\n\dot\ox(s^x_j)-\dot\ox(t_{j+1})\en\)\]
^2+2\Tilde\nu 4^{-k}\\
&\le 2h_k{\rm var}^2\(\dot\ox(\cdot);[0,T]\)+2\Tilde\nu4^{-k}\le
2h_k\mu^2+2\Tilde\nu 4^{-k}.
\end{aligned}
$$
This finally brings us to the estimate
$$
\int^T_0\n\dot x^k(t)-\dot\ox(t)\en^2dt\le
6(M^x_2)^2m^2\dd_k^2\Tilde\nu+6L\Tilde\nu\delta^2_k+12\Bar M\(h_k\mu+\Tilde\nu
2^{-k}+\sqrt{T\mu_k}\)+4h_k\mu^2+4\Tilde\nu4^{-k},
$$
which justifies the $L^2$-strong convergence of $\nn\dot x^k(\cdot)\hnn$ to
$\dot{\ox}(\cdot)$ in the norm topology as claimed at Step~3.\\[1ex]
{\bf Step~4:} {\em Verifying the convergence of $b^k(\cdot)$ to $\ob(\cdot)$ in
$W^{1,2}([0,T];\R^m)$.} This is the last step in the proof of the theorem.
Picking any $t\in(t_j,t_{j+1}]$, we have, using \eqref{e:bjk} and \eqref{e:bk},
$$
\begin{aligned}
\left|b^k_i(t)-\ob_i(t)\right|&=\left|b^k_{ij}+\dfrac{t-t_j}{h^k_j}\(b^k_{i,
j+1}-b^k_{ij}\)-\ob(t)\right|\\
&\le\left|\ob_i(t_j)-\ob_i(t)\right|+\left|\la
a^k_i(t_j),x^k_j\ra-\la\oa_i(t_j),\ox(t_j)\ra\right|+\left|b^k_{i,j+1}-b^k_{ij}
\right|.
\end{aligned}
$$
Since $\ob(\cdot)$ is uniformly continuous on $[0,T]$, for any $\ve>0$ there
exists $\delta>0$ ensuring that
$$
\max\left\{|t-s|,h_k\right\}<\delta\Longrightarrow\n\ob(t)-\ob(s)\en\le\ve,
$$
which implies that $\left|\ob_i(t_j)-\ob_i(t)\right|\le\ve$. Furthermore, it
follows from \eqref{e:a-est} and \eqref{e:x-est1} that
$$
\begin{aligned}
\left|\la a^k_i(t_j),x^k_j\ra-\la\oa_i(t_j),\ox(t_j)\ra\right|&=\left|\la
a^k_i(t_j)-\oa_i(t_j),x^k_j\ra+\la\oa_i(t_j),x^k_j-\ox(t_j)\ra\right|\\
&\le\n a^k_i(t_j)-\oa_i(t_j)\en\cdot\n x^k_j\en+\n\oa_i(t_j)\en\n
x^k_j-\ox(t_j)\en\\
&\le M_1\dd_k+\vartheta_k,
\end{aligned}
$$
where $M_1>0$ is chosen so that $\n x^k_j\en\le M_1$ for all $j=0,\ldots,k-1$, $\delta_k$
was defined in \eqref{e:def-deltak}, and $\vartheta_k$ was defined in \eqref{e:x-est1}.
Consequently we have
$$
\begin{aligned}
\left|b^k_{i,j+1}-b^k_{ij}\right|&=\left|\ob_i(t_{j+1})-\ob_i(t_j)+\la
a^k_i(t_{j+1}),x^k_{j+1}\ra-\la\oa_i(t_{j+1}),\ox(t_{j+1})\ra-\la
a^k_i(t_j),x^k_j\ra+\la\oa_i(t_j),\ox(t_j)\ra\right|\\
&\le\left|\ob_i(t_{j+1})-\ob_i(t_j)\right|+\left|\la
a^k_i(t_{j+1}),x^k_{j+1}\ra-\la\oa_i(t_{j+1}),\ox(t_{j+1})\ra\right|+\left|\la
a^k_i(t_j),x^k_j\ra-\la\oa_i(t_j),\ox(t_j)\ra\right|\\
&\le\ve+2\(M_1\dd_k+\vartheta_k\),
\end{aligned}
$$
which justifies the fulfillment of the claimed estimate
$$
\left|b^k_i(t)-\ob_i(t)\right|\le 2\ve+3\(M_1\dd_k+\vartheta_k\)
$$
and thus justifies the uniform convergence of $\nn b^k(\cdot)\hnn$ to
$\ob(\cdot)$ on $[0,T]$, thanks to \eqref{e:def-deltak}, \eqref{mu-k}, and \eqref{e:x-est1}.\vspace*{-0.05in}

It remains to prove the $L^2$-strong convergence of $\dot b^k(\cdot)$ to
$\dot\ob(\cdot)$ on $[0,T]$. For any $t\in[t_j,t_{j+1})$ we get
\begin{equation*}
\begin{aligned}
\left|\dot
b_i^k(t)-\dot\ob_i(t)\right|&=\left|\dfrac{b^k_{i,j+1}-b^k_{ij}}{h^k_j}
-\dot\ob_i(t)\right|\le
\left|\dfrac{\ob_i(t_{j+1})-\ob_i(t_j)}{h^k_j}-\dot\ob_i(t)\right|\\
&+\left|\la\dfrac{a^k_i(t_{j+1})-a^k_i(t_j)}{h^k_j}-\dfrac{\oa_i(t_{j+1}
)-\oa_i(t_j)}{h^k_j},x^k_{j+1}\ra\right|+\left|
\la\dfrac{\oa_i(t_{j+1})-\oa_i(t_j)}{h^k_j},x^k_{j+1}-\ox(t_{j+1})\ra\right|\\
&+\left|\la
a^k_i(t_j),\dfrac{x^k_{j+1}-x^k_j}{h^k_j}-\dfrac{\ox(t_{j+1})-\ox(t_j)}{h^k_j}
\ra\right|+\left|\la
a^k_i(t_j)-\oa_i(t_j),\dfrac{\ox(t_{j+1})-\ox(t_j)}{h^k_j}\ra\right|\\
&\le\left|\dfrac{\ob_i(t_{j+1})-\ob_i(t_j)}{h^k_j}-\dot\ob_i(t)\right|+M_1\n\dot
a^k_i(t)-\dfrac{\oa_i(t_{j+1})-\oa_i(t_j)}{h^k_j}\en\\
&+\vartheta_k\n\dfrac{\oa_i(t_{j+1})-\oa_i(t_j)}{h^k_j}\en+\(1+\dd_k\)\n\dot
x^k(t)-\dfrac{\ox(t_{j+1})-\ox(t_j)}{h^k_j}\en+\dd_k
\n\dfrac{\ox(t_{j+1})-\ox(t_j)}{h_k}\en
\end{aligned}
\end{equation*}
due to \eqref{e:bk}, \eqref{e:a-dc}, \eqref{e:a-est}, and \eqref{e:x-est1}.
Since $\oal(\cdot)$ is a BV function, it follows that
$$
\n\oal_i(s)\en\le M_2:=\frac{1}{2}\[\n\oal_{i0}\en+\n\oal_i(T)\en+{\rm
var}\big(\oal_i(\cdot);[0,T]\big)\]\;\;\mbox{for all}\;\;s\in[0,T],
$$
and therefore
$\n\dfrac{\oa_i(t_{j+1})-\oa_i(t_j)}{h_k}\en=\dfrac{1}{h_k}\n\disp\int^{t_{j+1}}
_{t_j}\oal_i(s)ds\en\le M_2$. Arguing in the same way for the BV function
$\dot\ox(\cdot)$ shows that $\n\dfrac{\ox(t_{j+1})-\ox(t_j)}{h_k}\en\le M_3$
with some constant $M_3>0$.\\
Next we estimate the quantities
$\left|\dfrac{\ob_i(t_{j+1})-\ob_i(t_j)}{h^k_j}-\dot\ob_i(t)\right|$ and $\n\dot
x^k(t)-\dfrac{\ox(t_{j+1})-\ox(t_j)}{h^k_j}\en$. Observe that
\begin{equation}\label{e:b-est1}
\begin{aligned}
\left|\dfrac{\ob_i(t_{j+1})-\ob_i(t_j)}{h^k_j}-\dot\ob_i(t)\right|&\le\dfrac{1}{
h^k_j}\int^{t_{j+1}}_{t_j}\left|\obe_i(s)-\obe_i(t_j)\right|ds+
\left|\obe_i(t_j)-\obe_i(t)\right|\\
&=\dfrac{1}{h^k_j}\int^{t_{j+1}}_{t_j}f^b_{ij}(s)ds+f^b_{ij}(t)\le
2f^b_{ij}(s^b_{ij})+2^{-k+1}\\
&\le
2\[\left|\obe_i(s^b_{ij})-\obe_i(t_j)\right|+\left|\obe_i(t_{j+1})-\obe_i(s^b_{
ij})\right|\]+2^{-k+1},
\end{aligned}
\end{equation}
which allows us while arguing as above to get the estimates
\begin{equation*}
\begin{aligned}
\n\dot a_i^k(t)-\dfrac{\oa_i(t_{j+1})-\oa_i(t_j)}{h^k_j}\en &\le\n\dot
a_i^k(t)-\dot\oa_i(t)\en+\n\dot\oa_i(t)-\dfrac{\oa_i(t_{j+1})-\oa_i(t_j)}{h^k_j}
\en\\
&\le\n\dot
a_i^k(t)-\dot\oa_i(t)\en+2\big[\n\oal_i(s^a_{ij})-\oal_i(t_j)\en+\n\oal_i(t_{j+1
})-\oal_i(s^a_{ij})\en\big]+2^{-k+1},
\end{aligned}
\end{equation*}
\begin{equation}\label{e:x-est2}
\begin{aligned}
\n\dot x^k(t)-\dfrac{\ox(t_{j+1})-\ox(t_j)}{h^k_j}\en &\le\n\dot
x^k(t)-\dot\ox(t)\en+\n\dot\ox(t)-\dfrac{\ox(t_{j+1})-\ox(t_j)}{h^k_j} \en\\
&\le \n\dot
x^k(t)-\dot\ox(t)\en+2\big[\n\dot\ox(s^x_{j})-\dot\ox(t_j)\en+\n\dot\ox(t_{j+1}
)-\dot\ox(s^x_{j})\en\big]+2^{-k+1}.
\end{aligned}
\end{equation}
It then follows by combining all the estimates in
\eqref{e:b-est1}--\eqref{e:x-est2} that
$$
\begin{aligned}
\int^T_0&\left|\dot
b_i^k(t)-\dot\ob_i(t)\right|^2dt=\sum^{\nu(k)-1}_{j=0}\int^{t_{j+1}}_{t_j}
\left|\dfrac{b^k_{i,j+1}-b^k_{ij}}{h^k_j}-\dot\ob_i(t)\right|^2dt\\
&\le\sum^{\nu(k)-1}_{j=0}\int^{t_{j+1}}_{t_j}\big[1+M_1^2+(1+\dd_k)^2+1\big]
\bigg[\left|\dfrac{\ob_i(t_{j+1})-\ob_i(t_j)}{h^k_j}- \dot\ob_i(t)\right|^2+
\n\dot a^k_i(t)-\dfrac{\oa_i(t_{j+1})-\oa_i(t_j)}{h^k_j}\en^2\\&+
\vartheta_k^2\n\dfrac{\oa_i(t_{j+1})-\oa_i(t_j)}{h^k_j}\en^2+\n\dot
x^k(t)-\dfrac{\ox(t_{j+1})-\ox(t_j)}{h^k_j}\en^2+\dd^2_k\n\dfrac{\ox(t_{j+1}
)-\ox(t_j)}{h^k_j}\en^2\bigg]dt\\
&\le\big[2+M_1^2+\(1+\dd_k\)^2\big]\bigg[\sum^{\nu(k)-1}_{j=0}\int^{t_{j+1}}_{
t_j}\left|\dfrac{\ob_i(t_{j+1})-\ob_i(t_j)}{h^k_j} -
\dot\ob_i(t)\right|^2dt+\sum^{\nu(k)-1}_{j=0}\int^{t_{j+1}}_{t_j}\n\dot
a^k_i(t)-\dfrac{\oa_i(t_{j+1})-\oa_i(t_j)}{h^k_j}\en^2dt\\
&+\sum^{\nu(k)-1}_{j=0}\int^{t_{j+1}}_{t_j}\vartheta_k^2M_2^2dt+\sum^{\nu(k)-1}_
{j=0}\int^{t_{j+1}}_{t_j}\n\dot
x^k(t)-\dfrac{\ox(t_{j+1})-\ox(t_j)}{h^k_j}\en^2dt+\dd^2_k\sum^{\nu(k)-1}_{j=0}
\int^{t_{j+1}}_{t_j}\n\dfrac{\ox(t_{j+1})-\ox(t_j)}{h^k_j}
\en^2dt\bigg]\\&\le\big[2+M_1^2+\(1+\dd_k\)^2\big]\bigg[4h_k\sum^{\nu(k)-1}_{j=0
}\[\left|\obe_i(s^b_{ij})-\obe_i(t_j)\right|+\left|\obe_i(t_{j+1})-\obe_i(s^b_{
ij}) \right|+2^{-k}\]^2+\sum^{\nu(k)-1}_{j=0}\int^{t_{j+1}}_{t_j}\n\dot
a_i^k(t)-\dot\oa_i(t)\en^2dt\\
&+4h_k\sum^{\nu(k)-1}_{j=0}\[\n\oal_i(s^a_{ij})-\oal_i(t_j)\en+\n\oal_i(t_{j+1}
)-\oal_i(s^a_{ij})\en+2^{-k}\]^2
+\vartheta_k^2M^2_2\Tilde\nu+\sum^{\nu(k)-1}_{j=0}\int^{t_{j+1}}_{t_j}\n\dot
x^k(t)-\dot\ox(t)\en^2dt\\
&+4h_k\sum^{\nu(k)-1}_{j=0}\[\n\dot\ox(s^x_{j})-\dot\ox(t_j)\en+\n\dot\ox(t_{j+1
})-\dot\ox(s^x_{j})\en\]^2+\dd^2_k
\sum^{\nu(k)-1}_{j=0}\int^{t_{j+1}}_{t_j}M_3^2dt\bigg].
\end{aligned}
$$
Finally, we arrive at the relationships
$$
\begin{aligned}
\int^T_0&\left|\dot b_i^k(t)-\dot\ob_i(t)\right|^2dt
\le\big[2+M_1^2+\(1+\dd_k\)^2\big]\bigg\{8h_k\sum^{\nu(k)-1}_{j=0}\[
\left|\obe_i(s^b_{ij})-\obe_i(t_j)\right|+\left|\obe_i(t_{j+1})-
\obe_i(s^b_{ij})\right|\]^2\\&+4^{-k}h_k\nu(k)+\int^{T}_{0}\n\dot
a_i^k(t)-\dot\oa_i(t)\en^2dt+8h_k\sum^{\nu(k)-1}_{j=0}\big[\n\oal_i(s^a_{ij}
)-\oal_i(t_j)\en+\n\oal_i(t_{j+1})-\oal_i(s^a_{ij})\en\big]^2\\
&+4^{-k}h_k\nu(k)+\vartheta_k^2M^2_2\Tilde\nu+\int^T_0\n\dot
x^k(t)-\dot\ox(t)\en^2dt\\
&+8h_k\sum^{\nu(k)-1}_{j=0}\big[\n\dot\ox(s^x_{j})-\dot\ox(t_j)\en+\n\dot\ox(t_{
j+1})-\dot\ox(s^x_{j})\en\big]^2+4^{-k}h_k\nu(k)+\dd^2_kM^2_3
\Tilde\nu\bigg\}\\
&\le\big[2+M_1^2+\(1+\dd_k\)^2\big]\bigg[8h_k\(\var^2(\obe(\cdot);[0,T]
)+\var^2(\oal(\cdot);[0,T])+\var^2(\dot\ox(\cdot);[0,T])\)+
\dfrac{3}{4^k}
\Tilde\nu\\&+\vartheta_k^2M^2_2\Tilde\nu+\dd^2_kM^2_3\Tilde\nu+\int^{T}_{0}
\n\dot a_i^k(t)-\dot\oa_i(t)\en^2dt+\int^T_0\n\dot
x^k(t)-\dot\ox(t)\en^2dt\bigg]\\&\le\big[2+M_1^2+\(1+\dd_k\)^2\big]\[
24h_k\mu^2+\dfrac{3}{4^k}\Tilde\nu+\vartheta_k^2M^2_2\Tilde\nu+\dd^2_kM^2_3
\Tilde\nu+\mu_k+\int^T_0\n\dot x^k(t)-\dot\ox(t)\en^2dt\],
\end{aligned}
$$
which ensures the convergence of the sequence $\{\dot b^k(\cdot)\}$ to
$\dot\ob(\cdot)$ strongly in $L^2([0,T];\R^m)$ as claimed in Step~4. This
therefore completes the proof of the theorem. $\h$

As we see, the entire proof of the theorem is technically involved. It occurs
nevertheless that the most important and challenging task is the construction of
a sequence of piecewise linear functions $x^k(\cdot)$, which are feasible to the
discrete differential inclusion \eqref{e:a-dc}. The main point is in
approximating the continuous velocity
$\dot\ox(t_j)\in-F(\ox(t_j),\oa(t_j),\ob(t_j),\ou(t_j))$ by its discrete
counterpart $v^k_j\in -F(x^k(t_j),a^k(t_j),b^k(t_j),u^k(t_j))$, where the
velocity mapping $F$ is discontinuous. Using the construction of $v^k_j$  in
\eqref{e:dv} ensures that the distance between $\dot\ox(t_j)$ and $v^k_j$
converges to $0$ as $k\to\infty$, which is the key.\vspace{-0.2in}

\section{Discrete Approximation for Relaxed Local Minimizers}\label{disc1}
\setcounter{equation}{0}\vspace*{-0.1in}

The discrete approximation procedure and results developed in the previous
section do not require any relaxation stability and do not concern optimal
versus feasible solutions. The discrete approximation construction and the main
result of this section address {\em relaxed local minimizers} of the sweeping
optimal control problem $(P)$.

Let $\(\ox(\cdot),\oa(\cdot),\ob(\cdot),\ou(\cdot)\)$ be a given r.i.l.m., and
let $\Delta_k$ be the discrete mesh defined in \eqref{e:DP}. For all $k\in\N$ we
construct a sequence of approximating problems $(P_k)$ as follows:
\begin{equation}\label{(P_k)}
\begin{aligned}
&\textrm{minimize
}J_k[x^k,a^k,b^k,u^k]:=\vph(x^k_{\nu(k)})+\sum^{\nu(k)-1}_{j=0}
h^k_j\ell\bigg(t^k_j,x^k_j,a^k_j,b^k_j,u^k_j,\dfrac{x^k_{j+1}-x^k_i}{h^k_j},
\dfrac{a^k_{j+1}-a^k_j}{h^k_j},\dfrac{b^k_{j+1}-b^k_j}{h^k_j}\bigg)\\
&+\frac{1}{2}\sum_{j=0}^{\nu(k)-1}\int_{t^k_j}^{t^k_{j+1}}\n\(\dfrac{x^k_{j+1}
-x^k_j}{h^k_j},\dfrac{a^k_{j+1}-a^k_j}{h^k_j},
\dfrac{b^k_{j+1}-b^k_j}{h^k_j},u^k_j\)-\(\dot{\ox}(t),\dot{\oa}(t),\dot{\ob}(t),
\ou(t)\)\en^2dt
\end{aligned}
\end{equation}
over discrete quadruples $(x^k,a^k,b^k,u^k)$ represented by
\begin{eqnarray*}
(x^k,a^k,b^k,u^k):=(x^k_0,x^k_1,\ldots,x^k_{\nu(k)},a^k_0,a^k_1,\ldots,a^k_{
\nu(k)},b^k_0,b^k_1,\ldots,b^k_{\nu(k)},u^k_0,u^k_1,\ldots,
u^k_{\nu(k)-1})
\end{eqnarray*}
subject to the geometric and functional constraints given by
\begin{equation}
x^k_{j+1}\in x^k_j-h^k_jF(x^k_j,a^k_j,b^k_j,u^k_j),\;j=0,\ldots,\nu(k)-1,
\end{equation}
\begin{equation}
\la a^{k}_{i\nu(k)},x^k_{\nu(k)}\ra\le b^{k}_{i\nu(k)},\;i=1,\ldots,m,
\end{equation}
\begin{equation}\label{ini}
x^k_0=x_0\in C(0),\;a^k_0=a_0,\;b^k_0=b_0,\;u^k_0=\ou(0),
\end{equation}
\begin{equation}\label{ic1}
\sum_{j=0}^{\nu(k)-1}\int_{t^k_j}^{t^k_{j+1}}\n\(x^k_j,a^k_j,b^k_j,
u^k_j\)-\(\ox(t),\oa(t),\ob(t),\ou(t)\)\en^2dt\le\dfrac{\ve}{2},
\end{equation}
\begin{equation}\label{ic2}
\sum_{j=0}^{\nu(k)-1}\int_{t^k_j}^{t^k_{j+1}}\n\(\dfrac{x^k_{j+1}-x^k_j}{h^k_j},
\dfrac{a^k_{j+1}-a^k_j}{h^k_j},\dfrac{b^k_{j+1}-b^k_j}{h^k_j}\)-
\(\dot{\ox}(t),\dot{\oa}(t),\dot{\ob}(t)\)\en^2dt\le\dfrac{\ve}{2},
\end{equation}
\begin{equation}
u^k_j\in U,\;j=0,\ldots,\nu(k)-1,
\end{equation}
\begin{equation}\label{constr_a}
1-\dd_k\le\|a^k_{ij}\|\le 1+\dd_k,\;i=1,\ldots,m,\;j=0,\ldots,\nu(k),
\end{equation}
where $\ve>0$ is taken from Definition~\ref{Def2} of the relaxed intermediate
local minimizer $\(\ox(\cdot),\oa(\cdot),\ob(\cdot),\ou(\cdot)\)$, where $F$ is
defined in \eqref{F}, and where the perturbation sequence $\dd_k\dn 0$ as
$k\to\infty$ is constructed in the proof of Theorem~\ref{Th2} for the given
quadruple $\(\ox(\cdot),\oa(\cdot),\ob(\cdot),\ou(\cdot)\)$.

To proceed further, first we need to make sure that for each $k\in\N$
sufficiently large the discrete control problem $(P_k)$ defined in
\eqref{(P_k)}--\eqref{constr_a} admits an optimal solution. It is verified in
the next proposition.\vspace*{-0.1in}

\begin{proposition}[\bf existence of optimal solutions to discrete sweeping
control problems]\label{Th4} Under the assumptions in Theorem~{\rm\ref{Th2}}
holding along the given r.i.l.m.\
$\(\ox(\cdot),\oa(\cdot),\ob(\cdot),\ou(\cdot)\)$, each problem $(P_{k})$ for
all sufficiently large $k\in\N$ admits an optimal solution.
\end{proposition}\vspace*{-0.12in}
{\bf Proof.} It follows from Theorem~\ref{Th2} that the set of feasible
solutions of problem $(P_k)$ is nonempty for all large $k$. We see in addition
that this set is bounded due to the constraint structures in $(P_k)$.
Furthermore, the cost function in $(P_k)$ is obviously lower semicontinuous for
each $t^k_j\in\Delta_k$ due to (H5). To apply the classical Weierstrass
existence theorem in $(P_k)$, it remains to ensure that the feasible set in this
problem is closed. But it is a direct consequence of the constraint structures
in $(P_k)$ due to the robustness (closed-graph) property of the normal cone
mapping \eqref{NC}. Thus we arrive at the claimed existence result. $\h$

Now we are ready to establish the desired theorem on the strong convergence of
optimal solutions for $(P_k)$ to the given r.i.l.m.\ of the original sweeping
control problem $(P)$.\vspace*{-0.1in}

\begin{theorem}[\bf strong convergence of discrete optimal solutions]\label{Th5}
Let $(\ox(\cdot),\oa(\cdot),\ob(\cdot),\ou(\cdot))$ be an r.i.l.m.\ for problem
$(P)$, and let all the assumptions of Proposition~{\rm\ref{Th4}} be satisfied
for this quadruple. Suppose in addition that the terminal cost $\vph$ is
continuous around $\ox(T)$, that the running cost $\ell$ is continuous at
$\big(t,\ox(t),\oa(t),\ob(t),\ou(t),$ $\dot\ox(t),\dot\oa(t),\dot\ob(t)\big)$
for a.e.\ $t\in[0,T]$, and that $\ell\(\cdot,x,a,b,u,\dot x,\dot a,\dot b\)$ is
uniformly majorized around $\(\ox(\cdot),\oa(\cdot),\ob(\cdot),\ou(\cdot)\)$ by
a summable function on $[0,T]$. Take any sequence of optimal solutions
$\(\ox^k(\cdot),\oa^k(\cdot),\ob^k(\cdot),\ou^k(\cdot)\)$ to the discrete
problems $(P_k)$  and extend it to the entire interval $[0,T]$ piecewise
linearly for $\(\ox^k(\cdot),\oa^k(\cdot),\ob^k(\cdot)\)$ and piecewise
constantly for $\ou^k(\cdot)$. Then the extended sequence
$\(\ox^k(\cdot),\oa^k(\cdot),\ob^k(\cdot),\ou^k(\
cdot)\)$ converges to $\(\ox(\cdot),\oa(\cdot),\ob(\cdot),\ou(\cdot)\)$  as
$k\to\infty$ in the norm topology of $W^{1,2}([0,T];\R^n)\times
W^{1,2}([0,T];\R^{mn})\times W^{1,2}([0,T];\R^m )\times L^2([0,T];\R^d)$.
\end{theorem}\vspace*{-0.12in}
{\bf Proof}. Picking any sequence
$\(\ox^k(\cdot),\oa^k(\cdot),\ob^k(\cdot),\ou^k(\cdot)\)$ of extended optimal
solutions to $(P_k)$, we claim that
\begin{equation}\label{gamma}
\lim_{k\to\infty}\int^T_0\n\(\dot\ox^k(t),\dot\oa^k(t),\dot\ob^k(t),
\ou^k(t)\)-\(\dot\ox(t),\dot\oa(t),\dot\ob(t),\ou(t)\)\en^2dt=0,
\end{equation}
which clearly ensures the convergence of the quadruples
$\(\ox^k(\cdot),\oa^k(\cdot),\ob^k(\cdot),\ou^k(\cdot)\)$ to
$\(\ox(\cdot),\oa(\cdot),\ob(\cdot),\ou(\cdot)\)$ in the norm topology of
$W^{1,2}([0,T];\R^{n+mn+m})\times L^2([0,T];\R^d)$. To proceed, assume on the
contrary that the limit in \eqref{gamma}, along a subsequence (without
relabeling), equals to some $\gamma>0$. Then it follows from the weak
compactness of the unit ball in $L^2([0,T];\R^{n+mn+m+d})$ that there exist
functions $(v^{x}(\cdot),v^{a}(\cdot),v^{b}(\cdot),\tu(\cdot))\in
L^2([0,T];\R^{n+mn+m+d})$ for which the quadruples
$(\dot\ox^k(\cdot),\dot\oa^k(\cdot),\dot\ob^k(\cdot),\ou^k(\cdot))$ converges
weakly to $\(v^{x}(\cdot),v^{a}(\cdot),v^{b}(\cdot),\tu(\cdot)\)$ in the
corresponding spaces. Recall that Mazur's weak closure theorem and basic real
analysis yield the existence of sequences of convex combinations of these
quadruples that converge to
$(v^{x}(\cdot),v^{a}(\cdot),v^{b}(\cdot),\tu(\cdot))$ in the $L^2$-norm topology
with
their subsequences (no relabeling) converging to
$(v^{x}(t),v^{a}(t),v^{b}(t),\tu(t))$ for a.e.\ $t\in[0,T]$. Define further the
triple $(\tx(\cdot),\ta(\cdot),\tb(\cdot))\in W^{1,2}([0,T];\R^{n+mn+m})$ by
\begin{equation*}
\(\tx(t),\ta(t),\tb(t)\):=\(x_0,a_0,b_0\)+\int^t_0\(v^x(s),v^a(s),
v^b(s)\)ds\;\mbox{ for all }\;t\in[0,T].
\end{equation*}
Then $(\dot\tx(t),\dot\ta(t),\dot\tb(t))=(v^x(t),v^a(t),v^b(t))$ for a.e.\
$t\in[0,T]$, which ensures the convergence of
$(\ox^k(\cdot),\oa^k(\cdot),\ob^k(\cdot))$ to
$(\tx(\cdot),\ta(\cdot),\tb(\cdot))$ in the norm topology of
$W^{1,2}([0,T];\R^{n+mn+m})$. Observe that $\tu(t)\in\co U$ for a.e.\
$t\in[0,T]$ and that the limiting triple $(\tx(\cdot),\ta(\cdot),\tb(\cdot))$
satisfies the differential inclusion \eqref{conv} with $C(t)=\Tilde
C(t):=\nn x\in \R^n|\;\la\ta_i(t),x\ra\le\tb_i(t),\; i=1,\ldots,m \hnn$ for all
$t\in[0,T]$. Taking into account the convexity of the norm function and hence
its lower semicontinuity in the $L^2$-weak topology, we get by passing to the
limit in \eqref{ic1} and \eqref{ic2}, respectively, that
$$
\begin{aligned}
\int^T_0&\n\big(\tx(t),\ta(t),\tb(t),\tu(t)\big)-\big(\ox(t),\oa(t),\ob(t),
\ou(t)\big)\en^2dt\\
&\le\liminf_{k\to\infty}\sum_{j=0}^{\nu(k)-1}\int_{t^k_j}^{t^k_{j+1}}
\n\big(x^k_j,a^k_j,b^k_j,u^k_j\big)-\big(\ox(t),\oa(t),\ob(t),\ou(t)\big)
\en^2dt\le\dfrac{\ve}{2},
\end{aligned}
$$
$$
\begin{aligned}
\int^T_0&\n\(\dot\tx(t),\dot\ta(t),\dot\tb(t)\)-\(\dot\ox(t),\dot\oa(t),
\dot\ob(t)\)\en^2dt\\
&\le\liminf_{k\to\infty}\sum_{j=0}^{\nu(k)-1}\int_{t^k_j}^{t^k_{j+1}}\n\(\dfrac{
x^k_{j+1}-x^k_j}{h^k_j},\dfrac{a^k_{j+1}-a^k_j}{h^k_j},
\dfrac{b^k_{j+1}-b^k_j}{h^k_j}\)-\(\dot{\ox}(t),\dot{\oa}(t),\dot{\ob}
(t)\)\en^2dt\le\dfrac{\ve}{2}
\end{aligned}
$$
This implies that the limiting quadruple
$(\tx(\cdot),\ta(\cdot),\tb(\cdot),\tu(\cdot))$ belongs to the given
$\ve$-neighborhood of the r.i.l.m.\
$(\ox(\cdot),\oa(\cdot),\ob(\cdot),\ou(\cdot))$ in the space
$W^{1,2}([0,T];\R^{n+mn+m})\times L^2([0,T];\R^d)$. It is clear furthermore that
$\ta(\cdot)$ satisfies the pointwise constraint \eqref{h1.5}. Applying now
Theorem~\ref{Th2} to the r.i.l.m.\
$\(\ox(\cdot),\oa(\cdot),\ob(\cdot),\ou(\cdot)\)$ gives us a sequence
$(x^k(\cdot),a^k(\cdot),b^k(\cdot),u^k(\cdot))$ of the extended feasible
solutions to $(P_k)$ such that $x^k(\cdot),a^k(\cdot),b^k(\cdot)$ and
$u^k(\cdot)$ strongly approximate $\ox(\cdot),\oa(\cdot),\ob(\cdot)$ and
$\ou(\cdot)$ in $W^{1,2}([0,T];\R^{n+mn+m})$ and $L^2([0,T];\R^d)$ respectively.
It then follows from the imposed convexity of $\hat\ell_F$ and the optimality of
$(\ox^k(\cdot),\oa^k(\cdot),\ob^k(\cdot),\ou^k(\cdot))$ to $(P_k)$ that
\begin{equation}\label{hat}
\begin{aligned}
\hat J&\[\tx,\ta,\tb,\tu\]+\dfrac{\gg}{2}
=\vph\big(\tx(T)\big)+\int^T_0\hat\ell_F\(t,\tx(t),\ta(t),\tb(t),\tu(t),
\dot\tx(t),\dot\ta(t),\dot\tb(t)\)dt+\dfrac{\gg}{2}\\
&\le\liminf_{k\to\infty}\[\vph\(\ox^k_{\nu(k)}\)+h_{k}\sum^{\nu(k)-1}_{j=0}
\ell\bigg(t^k_j,\ox^k_j,\oa^k_j,\ob^k_j,u^k_j,\dfrac{\ox^k_{j+1}
-\ox^k_i}{h^k_j},\dfrac{\oa^k_{j+1}-\oa^k_j}{h^k_j},\dfrac{\ob^k_{j+1}-\ob^k_j}{
h^k_j}\bigg)+\dfrac{\gg}{2}\]\\
&=\liminf_{k\to\infty}J_k\[\ox^k,\oa^k,\ob^k,\ou^k\]\le
\liminf_{k\to\infty}J_k\[x^k,a^k,b^k,u^k\],
\end{aligned}
\end{equation}
which ensures, in particular, that the quadruple $(\tx,\ta,\tb,\tu)$ is feasible
for the relaxed problem $(R)$. On the other hand, the strong
convergence of $(x^k(\cdot),a^k(\cdot),b^k(\cdot),u^k(\cdot))$ to
$\(\ox(\cdot),\oa(\cdot),\ob(\cdot),\ou(\cdot)\)$
in $W^{1,2}([0,T];\R^{n+mn+n})\times L^2([0,T];\R^d)$ from Theorem~\ref{Th2} and
the imposed continuity assumptions on $\vph$ and $\ell$ imply that
$J_k\[x^k,a^k,b^k,u^k\]\to J[\ox,\oa,\ob,\ou]$ as $k\to\infty$. Combining it
with \eqref{hat} tells us that
\begin{equation*}
\hat J\[\tx,\ta,\tb,\tu\]<\hat J\[\tx,\ta,\tb,\tu\]+\dfrac{\gg}{2}\le
J\[\ox,\oa,\ob,\ou\]=\hat J\[\ox,\oa,\ob,\ou\],
\end{equation*}
which clearly contradicts the fact that
$(\ox(\cdot),\oa(\cdot),\ob(\cdot),\ou(\cdot))$ is an r.i.l.m.\ for problem
$(P)$ and hence verifies the limiting condition \eqref{gamma}. This completes
the proof of the theorem. $\h$\vspace*{-0.2in}

\section{Generalized Differentiation and Second-Order Calculations}\label{tools}
\setcounter{equation}{0}\vspace*{-0.1in}

Having in hands the strong approximation results of Theorem~\ref{Th5}, our
subsequent goal is to derive necessary optimality conditions for the
discrete-time approximating problems $(P_k)$ that provide constructive
suboptimality conditions for the original sweeping control problem $(P)$.
Looking at problem $(P_k)$ for each fixed number $k\in\N$, we see that it is a
finite-dimensional optimization problem with various types of constraints. The
most important and challenging of these constraints, that are characteristic for
sweeping differential and finite-difference inclusions, are described by {\em
graphs} of {\em normal cone mappings}. Such sets are {\em nonconvex} regardless
of the convexity and/or smoothness of the given data of $(P)$. To deal with the
problems under consideration, we need to employ appropriate constructions of
{\em generalized differentiation} in variational analysis with paying the major
attention to {\em second-order} ones. This section briefly reviews the concepts
and results of
generalized differentiation used in what follows. We are mainly based on
\cite{m18}, while related first-order constructions can be also found in
\cite{rw}.

Recall that for a set-valued (in particular, single-valued) mapping
$S\colon\R^n\tto\R^m$ the symbol
\begin{equation}\label{out_lim}
\Limsup_{x\to\ox}S(x):=\big\{z\in\R^m\big|\;\exists\;\textrm{ sequences
}\;x_k\to\ox,\,z_k\to z\;\textrm{ such that }\;z_k\in S(x_k),\,k\in\N\big\}
\end{equation}
signifies the (Kuratowski-Painlev\'e) {\em outer limit} of $S$ at $\ox$. Given a
nonempty set $\O\subset\R^n$ locally closed around $\ox\in\O$, the (Mordukhovich
basic/limiting) {\em normal cone} to $\O$ at $\ox$ is defined via the outer
limit \eqref{out_lim} by
\begin{equation}\label{nor_con}
N(\ox;\O)=N_\O(\ox):=\Limsup_{x\to\ox}\big\{{\rm
cone}\big[x-\Pi(x;\O)\big]\big\},
\end{equation}
where $\Pi(\ox;\O)$ stands for the Euclidean projection of $\ox$ onto $\O$ and
is defined by
\begin{equation*}
\Pi(\ox;\O):=\big\{y\in\O\big|\;\|\ox-y\|=d(\ox;\O)\big\},
\end{equation*}
and where `cone' denotes the conic hull of a set. If $\O$ is convex, the
limiting normal cone \eqref{nor_con} reduces to the normal cone of convex
analysis \eqref{NC}, but in general this cone is nonconvex. Nevertheless, in
vast generality the normal cone \eqref{nor_con} as well as the associated
subdifferential and coderivative constructions enjoy comprehensive {\em calculus
rules} based on variational and extremal principles of variational analysis; see
\cite{m-book,m18,rw} for more details.\vspace*{-0.05in}

Given further a set-valued mapping $S\colon\R^n\tto\R^m$ whose graph
\begin{equation*}
\gph S:=\big\{(x,y)\in\R^n\times\R^m\big|\;y\in S(x)\big\}
\end{equation*}
is locally closed around $(\ox,\oy)$, the {\em coderivative} of $S$ at
$(\ox,\oy)$ is defined by
\begin{equation}\label{coderivative}
D^*S(\ox,\oy)(u):=\big\{v\in\R^n\big|\;(v,-u)\in N\big((\ox,\oy);\gph
S\big)\big\},\quad u\in\R^m.
\end{equation}
If $S\colon\R^n\to\R^m$ is single-valued and continuously differentiable (${\cal
C}^1$-smooth) around $\ox$, we have
\begin{equation*}
D^*S(\ox)(u)=\big\{\nabla S(\ox)^*u\big\}\;\textrm{ for all }\;u\in \R^m
\end{equation*}
via the adjoint/transposed Jacobian matrix $\nabla S(\ox)^*$, where $\oy=S(\ox)$
is omitted.\vspace*{-0.05in}

For an extended-real-valued l.s.c.\ function $\phi\colon\R^n\to\R$ with the
domain and epigraph defined by
\begin{equation*}
\dom\phi:=\big\{x\in\R^n\big|\;\phi(x)<\infty\big\}\;\textrm{ and
}\;\epi\phi:=\big\{(x,\al)\in\R^{n+1}\big|\;\al\ge\phi(x)\big\},
\end{equation*}
the {\em first-order subdifferential} of $\phi$ at $\ox\in\dom\phi$ is generated
geometrically via \eqref{nor_con} as
\begin{equation*}
\partial\phi(\ox):=\big\{v\in\R^n\big|\;(v,-1)\in
N\big((\ox,\phi(\ox));\epi\phi\big)\big\};
\end{equation*}
see \cite{m-book,m18,rw} for equivalent analytic representations. The {\em
second-order subdifferential}, or {\em generalized Hessian}, of $\phi$ at $\ox$
relative to $\ov\in\partial\phi(\ox)$ is the mapping
$\partial^2\phi(\ox,\ov)\colon\R^n\tto\R^n$ with the values
\begin{equation}\label{2nd}
\partial^2\phi(\ox,\ov)(u):=\big(D^*\partial\phi\big)(\ox,\ov)(u),\quad
u\in\R^n.
\end{equation}
If $\phi$ is a ${\cal C}^2$-smooth around $\ox$, then \eqref{2nd} with
$\ov=\nabla\phi(\ox)$ reduces to to the classical (symmetric) Hessian matrix:
\begin{equation*}
\partial^2\phi(\ox,\ov)(u)=\big\{\nabla^2\phi(x)u\big\}\;\mbox{ for all
}\;u\in\R^n.
\end{equation*}
Our main interest in this paper corresponds to the case where
$\phi(x)=\dd_\O(x)$ is the indicator function of a set that equals to $0$ for
$x\in\O$ and $\infty$ otherwise. In this case we have
$\partial\dd_\O(\ox)=N_\O(\ox)$ whenever $\ox\in\O$. The following result
presents evaluations of the coderivative \eqref{coderivative} of the normal cone
mapping
\begin{equation}\label{G}
G\colon\R^n\times\R^{mn}\times\R^m\tto\R^n\;\mbox{ with
}\;G(x,a,b):=N\big(x;C(a,b)\big)
\end{equation}
associated with the moving polyhedral set \eqref{C}. In fact, we get an
efficient upper estimate of the coderivative under PLICQ \eqref{PLICQ} and its
precise calculation under the corresponding LICQ. The proof of this result,
given in \cite[Lemmas~4.1 and 4.2]{chhm3}, is based on the second-order calculus
obtained in \cite{mo} and the seminal theorem by Robinson \cite{rob} on the
upper Lipschitzian stability of polyhedral multifunctions. To proceed, consider
the matrix $A:=[a_{ij}]$ as $i=1,\ldots,m$ and $j=1,\ldots,n$ with the vector
columns $a_i$, $i=1,\ldots,n$. Recall that the symbol $^\perp$ indicates the
orthogonal complement of a vector in the space in question.\vspace*{-0.1in}

\begin{lemma}[\bf coderivative evaluations of the normal cone
mapping]\label{Th6} Let $G$ be defined in \eqref{G} with $x\in C(a,b)$ for
$(x,a,b)\in\R^n\times\R^{mn}\times\R^m$, and let $v\in G(x,a,b)$. Suppose that
the active constraint vectors $\nn a_i|\;i\in I(x,a,b)\hnn$ are positively
linearly independent. Then we have the coderivative upper estimate
$$
D^*G(x,a,b,v)(w)\subset\bigcup\nn\(\begin{array}{c}
A^*q\\
\hline
p_1w+q_1x\\
\vdots\\
p_mw+q_mx\\
\hline
-q
\end{array}\)
\Bigg|\;p\in N_{\R^m_-}(Ax-b),\;A^*p=v,\;q\in D^*N_{\R^m_-}(Ax-b,p)(Aw)\hnn
$$
for $w\in\disp\bigcap_{\nn i|\;p_i>0\hnn}a^\perp_i\subset\R^n$ and
$D^*G(x,a,b,v)(w)=\emp$ otherwise.\\
If the active constraint vectors $\nn a_i|\;i\in I(x,a,b)\hnn$ are linearly
independent, then we have the precise formula
$$
D^*G(x,a,b,v)(w)=\bigcup_{q\in D^*N_{\R^m_-}(Ax-b,p)(Aw)}\(\begin{array}{c}
A^*q\\
\hline
p_1w+q_1x\\
\vdots\\
p_mw+q_mx\\
\hline
-q
\end{array}\)\;\;\mbox{for all}\;\;w\in\disp\bigcap_{\nn
i\big|\;p_i>0\hnn}a^\perp_i,
$$
where the vector $p\in N_{\R^m_-(Ax-b)}$ is uniquely determined by $A^*p=v$.
Furthermore, the coderivative of the normal cone mapping \eqref{G} generated by
the nonpositive orthant $\R^m_-$ is computed by
\begin{eqnarray}\label{e:cor-nc-o}
D^*N_{\R^m_-}(x,v)(w)=\left\{\begin{array}{ll}
\emp&\mbox{if }\,\exists\;i\;\mbox{ with }\;v_iw_i\not=0\\
\big\{\gg\big|\;\gg_i=0\;\forall\;i\in I_1(w),\;\gg_i\ge 0\;\forall\;i\in
I_2(w)\big\}&\mbox{otherwise}
\end{array}\right.,
\end{eqnarray}
whenever $(x,v)\in\gph N_{\R^m_-}$ with the index subsets in \eqref{e:cor-nc-o}
defined by
\begin{equation}\label{e:indexes}
I_1(w):=\big\{i\big|\;x_i<0\big\}\cup\big\{i\big|\;v_i=0,\;w_i<0\big\},\quad
I_2(w):=\big\{i\big|\;x_i=0,\;v_i=0,\;w_i>0\big\}.
\end{equation}
\end{lemma}

The following theorem, which is is strongly used in deriving necessary
optimality conditions in the next section, provides constructive evaluations of
the coderivative of the sweeping control mapping $F$ taken from \eqref{F}
entirely in terms of the given problem data. \vspace*{-0.1in}

\begin{theorem}[\bf coderivative evaluations of the sweeping control
mapping]\label{Th7} Consider the multifunction $F$ from \eqref{F} with the
polyhedral set $C$ defined in \eqref{C}, where the perturbation mapping $g(x,u)$
is ${\cal C}^1$-smooth around the reference points, and where $G$ is defined in
\eqref{G}. Suppose that the vectors $\nn a_i|\;i\in I(x,a,b)\hnn$ are positively
linearly independent for any triple $(x,a,b)\in\R^n\times\R^{mn}\times\R^m$.
Then for all such triples and all $(w,u)\in\R^m\times U$ with $w +g(x,u)\in
G(x,a,b)$ we have the coderivative upper estimate
\begin{equation}\label{c57}
\begin{aligned}
D^*F(x,a,b,u,w)(y)\subset\bigcup_{\substack{p\in
N_{\R^m_-}(Ax-b),\;A^*p=w+g(x,u)\\q\in D^*N_{\R^m_-}(Ax-b,p)(Ay)}}
\(\begin{array}{c}
A^*q-\nabla_x g(x,u)^*y\\
\hline
p_1y+q_1x\\
\vdots\\
p_my+q_mx\\
\hline
-q\\
-\nabla_u g(x,u)^*y
\end{array}\)
\end{aligned}
\end{equation}
for any $y\in\disp\bigcap_{\nn i|\;p_i>0\hnn}a^\perp_i$, where the vector
$q\in\R^m$ satisfies the conditions
\begin{equation}\label{q-c}
\left\{\begin{array}{ll}
q_i=0\;\mbox{ for all }\;i\;\mbox{ such that either }\;\la a_i,x\ra<b_i\;\mbox{
or }\;p_i=0,\;\mbox{ or}\;\la a_i,y\ra<0,\\
q_i\ge 0\;\mbox{ for all }\;i\;\mbox{ such that }\;\la
a_i,x\ra=b_i,\;p_i=0,\;\mbox{ and }\;\la a_i,y\ra>0.
\end{array}
\right.
\end{equation}
Furthermore, the equality holds in \eqref{c57} if the vectors $\nn a_i|\;i\in
I(x,a,b)\hnn$ are linearly independent in which case the vector $p\in
N_{\R^m_-(Ax-b)}$ is uniquely determined by $A^*p=w+g(x,u)$.
\end{theorem}\vspace*{-0.15in}
{\bf Proof.} Pick any $y\in\disp\bigcap_{\nn i|\;p_i>0\hnn}a^\perp_i$ and any
$z^*\in D^*F(x,a,b,u,w)(y)$. It follows from the coderivative sum rules of the
equality type given in \cite[Theorem~3.9]{m18} that
\begin{equation*}
z^*\in\(\begin{array}{c}
-\nabla_x g(x,u)^*y\\
\hline
0\\
\vdots\\
0\\
\hline
-\nabla_u g(x,u)^*y
\end{array}
\)y+\(\begin{array}{c|c}
I&0\\
\hline
0&0
\end{array}\)
D^*G\big(x,a,b,w+g(x,u)\big)(y).
\end{equation*}
Employing further Lemma~\ref{Th6} tells us that
$$
z^*\in\(\begin{array}{c}
-\nabla_x g(x,u)^*y\\
\hline
0\\
\vdots\\
0\\
\hline
-\nabla_u g(x,u)^*y
\end{array}
\)y+\(\begin{array}{c|c}
I&0\\
\hline
0&0
\end{array}\)
\(\begin{array}{c}
A^*q\\
\hline	
p_1y+q_1x\\
\vdots\\
p_my+q_mx\\
\hline
-q
\end{array}\)=\(\begin{array}{c}
A^*q-\nabla_x g(x,u)^*y\\
\hline
p_1y+q_1x\\
\vdots\\
p_my+q_mx\\
\hline
-q\\
-\nabla_u g(x,u)^*y
\end{array}\)
$$
for some $p\in N_{\R^m_-}(Ax-b)$ with $A^*p=w+g(x,u)$ and $ q\in
D^*N_{\R^m_-}(Ax-b,p)(Ay)$. Finally, conditions \eqref{q-c} for the vector $q$
follows from \eqref{e:cor-nc-o} and \eqref{e:indexes}. This completes the proof
of the theorem. $\h$\vspace*{-0.15in}

\section{Optimality Conditions via Discrete Approximations }\label{nec}
\setcounter{equation}{0}\vspace*{-0.1in}

This section is devoted to deriving necessary optimality conditions for each
discrete-time problem $(P_k)$ as $k\in\N$. As follows from Theorem~\ref{Th5},
the results obtained below give us suboptimality conditions for the selected
r.i.l.m.\ of the original sweeping optimal control problem $(P)$ provided that
the discretization index $k$ is sufficiently large.\vspace*{-0.02in}

We establish here two results in this direction. The first theorem provides
necessary optimality conditions to each problem $(P_k)$ defined in
Section~\ref{disc1} that are expressed in terms of the normal cone to the graph
of the velocity mapping $F$ from \eqref{F}, i.e., via the coderivative of this
mapping. The second theorem is the main result of this section. It derives
verifiable necessary conditions for the given r.i.l.m.\ of problem $(P)$
expressed entirely in terms of the initial data of the original sweeping control
problem along the strongly converging sequence of optimal solutions to the
discrete approximation problems $(P_k)$.

Let us start with the first result, which proof is based on the reduction of
$(P_k)$ to nonsmooth finite-dimensional mathematical programming with
increasingly many geometric constraints and employing calculus rules of
first-order generalized differentiation. As seen below, the proof of the main
result is largely based on second-order calculations. For convenience we use the
notation $\rep_m(x):=\(x,\ldots,x\)\in\R^{mn}$.\vspace*{-0.07in}

\begin{theorem}[\bf necessary conditions for discrete optimal
solutions]\label{Th8} Fix any $k\in\N$ and let
\begin{equation*}
(\ox^k,\oa^k,\ob^k,\ou^k)=(\ox_0^k,\ldots,\ox_{\nu(k)}^k,\oa_0^k,\ldots,\oa_{
\nu(k)}^k,\ob_0^k,\ldots,\ob_{\nu(k)}^k,\ou_0^k,\ldots,\ou_{\nu(k)-1}^k)
\end{equation*}
be an optimal solution to $(P_k)$ along which the general assumptions of
Theorem~{\rm\ref{Th7}} are fulfilled. Suppose in addition that the cost
functions $\vph$ and $\ell$ are locally Lipschitzian around the corresponding
components of the optimal solution. Then there exist a number $\lm^k\ge 0$ and
vectors
$\al^{1k}=\(\al_{0}^{1k},\ldots,\al_{\nu(k)}^{1k}\)\in\R^{(\nu(k)+1)m}_+$,
$\psi^k=(\psi^k_0,\ldots,\psi^k_{\nu(k)-1})\in\R^{\nu(k)d}$,
$\al^{2k}=\(\al_{0}^{2k},\ldots,\al_{\nu(k)}^{2k}\)\in\R^{(\nu(k)+1)m}_{-}$,
$\xi^k=(\xi_1^k,\ldots,\xi_{m}^k)\in\R^{m}_+$, and
$p_j^k=\(p_{j}^{xk},p_{j}^{ak},p_{j}^{bk}\)\in\R^{n+mn+m}$ as
$j=0,\ldots,\nu(k)$ satisfying the relationships:
\begin{equation}\label{e:5.8*}
\lm^k+\n\xi^k\en+\|\al^{1k}+\al^{2k}\|+\sum_{j=0}^{\nu(k)-1}\n
p_{j}^{xk}\en+\|p_{0}^{ak}\|+\|p_{0}^{bk}\|+\n\psi^k\en\ne 0,
\end{equation}
\begin{equation}\label{xi}
\xi_i^k\(\la\oa_{ik}^k,\ox_{k}^k\ra-\ob_{ik}^k\)=0,\;\;i=1,\ldots,m,
\end{equation}
\begin{equation}\label{al1}
\al_{ij}^{1k}\(\|\oa_{ij}^k\|-(1+\dd_k)\)=0,\;\;i=1,\ldots,m,\;\;j=0,\ldots,
\nu(k),
\end{equation}
\begin{equation}\label{al2}
\al_{ij}^{2k}\(\|\oa_{ij}^k\|-(1-\dd_k)\)=0,\;i=1,\ldots,m,\;\;j=0,\ldots,\nu(k)
,
\end{equation}
\begin{equation}\label{mutx}
-p_{\nu(k)}^{xk}=\lm^kv_{\nu(k)}^k+\sum_{i=1}^{m}\xi_i^k\oa_{i\nu(k)}
^k\in\lm^k\partial\vph(\ox_{\nu(k)}^k)+\sum_{i=1}^{m}\xi_i^k\oa_{i\nu(k)}^k,
\end{equation}
\begin{equation}\label{muta}
p_{\nu(k)}^{ak}=-2\[\al_{\nu(k)}^{1k}+\al_{\nu(k)}^{2k},\oa_{i\nu(k)}^k\]-\[
\xi^k,\rep_m(\ox_{\nu(k)}^k)\],
\end{equation}
\begin{equation}\label{mutb}
p_{\nu(k)}^{bk}=\xi^k,
\end{equation}
\begin{equation}\label{e:dac5}
p^{ak}_{j+1}=\lm^k\(v^{ak}_j+\dfrac{1}{h^k_j}\th^{Ak}_j\),\;\;p^{bk}_{j+1}
=\lm^k\(v^{bk}_j+\dfrac{1}{h^k_j}\th^{Bk}_j\),\;\;j=0,\ldots,\nu(k)-1,
\end{equation}
\begin{equation}\label{e:5.10*}
\begin{aligned}
&\Bigg(\frac{p_{j+1}^{xk}-p_{j}^{xk}}{h^k_j}-\lm^kw_{j}^{xk},\frac{p_{j+1}^{ak}
-p_{j}^{ak}}{h^k_j}-\lm^kw_{j}^{ak},
\frac{p_{j+1}^{bk}-p_{j}^{bk}}{h^k_j}-\lm^kw_{j}^{bk},-\frac{1}{h_k}\lm^k\th_{j}
^{uk}-\lm^kw_{j}^{uk},-p_{j+1}^{xk}\\
&+\lm^k\Big(v_{j}^{xk}+\dfrac{1}{h^k_j}\th_{j}^{Xk}\Big)\Bigg)\in\(0,\frac{2}{
h^k_j}
\[\al_{j}^{1k}+\al_{j}^{2k},\oa^k_j\],0,\frac{1}{h^k_j}\psi^k_j,0\)\\
&+N\(\(\ox_j^k,\oa_j^k,\ob_j^k,\ou_j^k,-\frac{\ox_{j+1}^k-\ox_j^k}{h^k_j}\);\gph
F\),\;\;j=0,\ldots,\nu(k)-1,
\end{aligned}
\end{equation}
\begin{equation}\label{psi}
\psi^k_j\in N\(u^k_j;U\),\;\;j=0\ldots,\nu(k)-1,
\end{equation}
where the quadruple $(\th_{j}^{uk},\th_{j}^{Xk},\th_{j}^{Ak},\th_{j}^{Bk})$ is
defined by
\begin{equation*}
\left(\int_{t_j^k}^{t_{j+1}^k}\(\ou_j^k-\ou(t)\)dt,\int_{t_j^k}^{t_{j+1}^k}
\(\frac{\ox_{j+1}^k-\ox_j^k}{h^k_j}-\dot{\ox}(t)\)dt,
\int_{t_j^k}^{t_{j+1}^k}\(\frac{\oa_{j+1}^k-\oa_j^k}{h^k_j}-\dot{\oa}(t)\)dt,
\int_{t_j^k}^{t_{j+1}^k}\(\frac{\ob_{j+1}^k-\ob_j^k}{h^k_j}
-\dot{\ob}(t)\)dt\right)
\end{equation*}
with the running cost subgradient collections
\begin{equation}\label{subcol}
\(w_{j}^{xk},w_{j}^{ak},w_{j}^{bk},w_{j}^{uk},v_{j}^{xk},v_{j}^{ak},v_{j}^{bk}\)
\in\partial\ell\(\ox_j^k,\oa_j^k,\ob_j^k,\ou_j^k,\frac{\ox_{j+1}^k-\ox_j^k}{
h^k_j},\frac{\oa_{j+1}^k-\oa_j^k}{h^k_j},\frac{\ob_{j+1}^k
-\ob_j^k}{h^k_j}\).
\end{equation}
\end{theorem}\vspace*{-0.1in}
{\bf Proof.} Take $\ve>0$ from the definition of the r.i.l.m.\
$(\ox(\cdot),\oa(\cdot),\ob(\cdot),\ou(\cdot))$ in problem $(P_k)$ for any fixed
$k\in\N$ and define the ``long" vector reflecting the set of feasible solutions
to each discrete-time problem $(P_k)$ by
$$
z:=\(x_0^k,\ldots,x_{\nu(k)}^k,a_0^k,\ldots,a_{\nu(k)}^k,b_0^k,\ldots,b_{\nu(k)}
^k,u_0^k,\ldots,u_{\nu(k)-1}^k,X_0^k,\ldots,X_{\nu(k)-1}^k,A_0^k,\ldots,A_{
\nu(k)-1}^k,B_0^k,\ldots,B_{\nu(k)-1}^k\)
$$
with the fixed starting point as in \eqref{ini}. It is clear that each problem
$(P_k)$ can be equivalently written as the nondynamic problem of mathematical
programming $(MP)$ with respect to vector $z$:
\begin{equation}\label{e:5.12*}
\begin{aligned}
&\textrm{minimize
}\phi_0(z):=\vph\big(x(T)\big)+\sum_{j=0}^{\nu(k)-1}h^k_j\ell(x_j^k,a_j^k,b_j^k,
u_j^k,X_j^k,A_j^k,B_j^k)\\
&+\frac{1}{2}\sum_{j=0}^{\nu(k)-1}\int_{t_j^k}^{t_{j+1}^k}\n\(X_j^k-\dot{\ox}(t)
,A_j^k-\dot{\oa}(t),B_j^k-\dot{\ob}(t),u_j^k-\ou(t)\)\en^2dt
\end{aligned}
\end{equation}
subject to finitely many equality, inequality, and geometric constraints
\begin{equation}\label{e:5.13*}
\kappa(z):=\sum_{j=0}^{\nu(k)-1}\int_{t^k_j}^{t^k_{j+1}}\n\(x^k_j,a^k_j,b^k_j,
u^k_j\)-\(\ox(t),\oa(t),\ob(t),\ou(t)\)\en^2dt\le\dfrac{\ve}{2},
\end{equation}
\begin{equation}\label{e:5.14*}
\phi(z):=\sum_{j=0}^{\nu(k)-1}\int_{t_j^k}^{t_{j+1}^k}\n\(X_j^k,A_j^k,B_j^k,
u_j^k\)-(\dot{\ox}(t),\dot{\oa}(t),\dot{\ob}(t),\ou(t))\en^2dt-\frac{\e}{2}\le
0,
\end{equation}
\begin{equation}
g^x_j(z):=x_{j+1}^k-x_j^k-h^k_jX_j^k=0,\;j=0,\ldots,\nu(k)-1,
\end{equation}
\begin{equation}
g^a_j(z):=a_{j+1}^k-a_j^k-h^k_jA_j^k=0,\;j=0,\ldots,\nu(k)-1,
\end{equation}
\begin{equation}
g^b_j(z):=b_{j+1}^k-b_j^k-h^k_jB_j^k=0,\;j=0,\ldots,\nu(k)-1,
\end{equation}
\begin{equation}
q_i(z):=\la a_{i\nu(k)}^k,x_{\nu(k)}^k\ra-b_{i\nu(k)}^k\leq 0,\;i=1,\ldots,m,
\end{equation}
\begin{equation}
l^1_{ij}(z):=\n a_{ij}^k\en^2-(1+\dd_k)^2\le
0,\;i=1,\ldots,m,\;j=0,\ldots,\nu(k),
\end{equation}
\begin{equation}
l^2_{ij}(z):=\n a_{ij}^k\en^2-(1-\dd_k)^2\ge
0,\;i=1,\ldots,m,\;j=0,\ldots,\nu(k),
\end{equation}
\begin{equation}\label{e:5.16*}
z\in\Xi_j:=\left\{z\big|\;-X_j^k\in
F(x_j^k,a_j^k,b_j^k,u_j^k)\right\},\;j=0,\ldots,\nu(k)-1,
\end{equation}
\begin{equation}\label{e:5.17*}
z\in\Xi_{\nu(k)}:=\left\{z\big|\;x_0^k\textrm{ is
fixed},\;\(a_0,b_0,u_0\)=\(\oa(0),\ob(0),\ou(0)\)\right\},
\end{equation}
\begin{equation}\label{uinU}
z\in\O_j=\big\{z\big|\;u^k_j\in U\},\;\;j=0,\ldots,\nu(k)-1.
\end{equation}
Next we apply the necessary conditions from \cite[Proposition~6.4(ii) and
Theorem~6.5(ii)]{m-book} to the optimal solution
\begin{eqnarray*}
\oz:=\(\ox_0^k,\ldots,\ox_{\nu(k)}^k,\oa_0^k,\ldots,\oa_{\nu(k)}^k,\ob_0^k,
\ldots,\ob_{\nu(k)}^k,\ou_0^k,\ldots,\ou_{\nu(k)-1}^k,\oX_0^k,\ldots,\oX_{
\nu(k)-1}^k,\oA_0^k,\ldots,\oA_{\nu(k)-1}^k,\oB_0^k,
\ldots,\oB_{\nu(k)-1}^k\)
\end{eqnarray*}
of problem $(MP)$ in \eqref{e:5.12*}--\eqref{uinU} corresponding to
the one for $(P_k)$ given in the theorem. It follows from Theorem~\ref{Th5} that
the inequality constraints in \eqref{e:5.13*} and \eqref{e:5.14*} are inactive
for large $k$, and so the corresponding multipliers do not appear in the
necessary optimality conditions. Thus we find dual elements $\lm^k\ge
0,\;\xi^k=(\xi_1^k,\ldots,\xi_m^k)\in\R^m_+,\;\al^{1k}=(\al_{0}^{1k},\ldots,\al_
{\nu(k)}^{1k})\in\R^{\nu(k)+1}_+,\;\al^{2k}=(\al_{0}^{2k},\ldots,\al_{\nu(k)}^{
2k})\in\R^{\nu(k)+1}_{-},\;p_j^k=
\(p_{j}^{xk},p_{j}^{ak},p_{j}^{bk}\)\in\R^{n+mn+m}$ for $j=1,\ldots,\nu(k)$, and
$z^*_j=\big(x^*_{0j},\ldots,x^*_{\nu(k)j},a^*_{0j},\ldots,a^*_{\nu(k)j},b^*_{0j}
,\ldots,b^*_{\nu(k)j},u^*_{0j},\ldots,u^*_{(\nu(k)-1)j},X^*_{0j},\ldots,\\X^*_{
(\nu(k)-1)j},A^*_{0j},\ldots,A^*_{(\nu(k)-1)j},
B^*_{0j},\ldots,B^*_{(\nu(k)-1)j}\big)$ for $j=0,\ldots,\nu(k)$, which are not
zero simultaneously, such the the following relationships are satisfied:
\begin{equation}
\label{69}
z^*_j\in\left\{
\begin{array}{llll}
N(\oz,\Xi_j)+N(\oz,\O_j)&\textrm{ if }&j\in\{0,\ldots,\nu(k)-1\}\\
N(\oz,\Xi_j)&\textrm{ if }& j=\nu(k)
\end{array}
\right.,
\end{equation}
\begin{equation}\label{70}
-z^*_0-\ldots-z^*_{\nu(k)}\in\lm^k\partial\phi_0(\oz)+\sum_{i=1}^{m}
\xi_i^k\nabla
q_i(\oz)+\sum_{j=0}^{\nu(k)}\sum^m_{i=1}\al_{ij}^{1k}\nabla
l^1_{ij}(\oz)+\sum_{j=0}^{\nu(k)}\sum^m_{i=1}\al_{ij}^{2k}\nabla
l^2_{ij}(\oz)+\sum_{j=0}^{\nu(k)-1}(\nabla g_j(\oz))^*p_{j+1}^{k},
\end{equation}
\begin{equation}\label{71+}
\xi_i^k q_i(\oz)=0,\;\;i=1,\ldots,m,
\end{equation}
\begin{equation}\label{71l1}
\al_{ij}^{1k}\(\n a_{ij}^k\en
-(1+\dd_k)\)=0,\;i=1,\ldots,m,\;\;j=0,\ldots,\nu(k),
\end{equation}
\begin{equation}\label{71l2}
\al_{ij}^{2k}\(\n a_{ij}^k\en -(1-\dd_k)\)=0,\;i=1,\ldots,m,\;j=0,\ldots,\nu(k).
\end{equation}
Note that the first line in $(\ref{69})$ comes from applying the normal cone
intersection rule from \cite[Theorem~2.16]{m18} to $\oz\in\O_j\cap\Xi_j$ for
$j=0,\ldots,\nu(k)-1$, where the qualification condition
\begin{equation}\label{qc}
N(\oz;\Xi_j)\cap\big(-N(\oz;\O_j)\big)=\{0\},\quad j=0,\ldots,\nu(k)-1,
\end{equation}
imposed therein is fulfilled. Indeed, for any vector $z^*_j\in
N(\oz;\Xi_j)\cap(-N(\oz;\O_j))$ we clearly have the inclusions
\begin{equation}\label{qc1}
(x^*_{jj},a^*_{jj},b^*_{jj},u^*_{jj},-X^*_{jj})\in
N\Big(\Big(\ox^k,\oa^k_j,\ob^k_j,\ou^k_j,-\frac{\ox^{k+1}_j-\ox^{k+1}_j}{h_k}
\Big);\gph F\Big),\quad-u^*_{ij}\in N(\ou^k_j;U),
\end{equation}
while the other components of $z^*_j$ are zero. It immediately follows from the
second inclusion in \eqref{qc1} that
\begin{equation*}
x^*_{jj}=0,\;a^*_{jj}=0,\;b^*_{jj}=0,\;\mbox{ and }\;X^*_{jj}=0.
\end{equation*}
Substituting this into the first inclusion in \eqref{qc1} and using the
coderivative definition \eqref{coderivative} give us
\begin{equation*}
(0,0,0,u^*_{jj})\in
D^*F\Big(\ox^k_j,\oa^k_j,\ob^k_j,\ou^k_j,-\frac{\ox^{k+1}_j-\ox^k_j}{h_k}
\Big)(0),\quad j=0,\ldots,\nu(k)-1.
\end{equation*}
Then we deduce directly from the coderivative estimate \eqref{c57} for the
velocity mapping $F$ in \eqref{F} under the imposed PLICQ that $u^*_{jj}=0$ for
all $j=0,\ldots,\nu(k)-1$. It shows that $z^*_j=0$ for such indices $j$, and
therefore the qualification condition \eqref{qc} is verified.

To proceed further, observe from the structure of the sets $\Xi_j$ and $\O_j$ in
\eqref{e:5.16*}--\eqref{uinU}, respectively, that the inclusions in \eqref{69}
are equivalent to
\begin{equation}\label{e:5.18*}
\left\{\begin{array}{ll}
(x^*_{jj},a^*_{jj},b^*_{jj},u^*_{jj}-\psi_{j}^k,-X^*_{jj})\in
N\Big(\Big(\ox_j^k,\oa_j^k,\ob_j^k,\ou_j^k,-\dfrac{\ox_{j+1}^k-\ox_j^k}{h^k_j}
\Big);\gph F\Big)&\textrm{ for }\;j=0,\ldots,\nu(k)-1,\\
(x^*_{jj},a^*_{jj},b^*_{jj},u^*_{jj},- X^*_{jj})\in
N\Big(\Big(\ox_j^k,\oa_j^k,\ob_j^k,\ou_j^k,-\dfrac{\ox_{j+1}^k-\ox_j^k}{h^k_j}
\Big);\gph F\Big)&\textrm{ for }\;j=\nu(k)
\end{array}
\right.
\end{equation}
with $\psi_{j}^k$ taken from \eqref{psi}, while the other components of $z^*_j$
equal to zero. Similarly we get that the vectors
$x^*_{0\nu(k)},a^*_{0\nu(k)},b^*_{0\nu(k)}$, and $u^*_{0\nu(k)}$ determined by
the normal cone to $\Xi_{\nu(k)}$ might be the only nonzero components of
$z^*_{\nu(k)}$. This readily yields the representation
\begin{eqnarray*}
\begin{aligned}
&-z^*_0-\ldots-z^*_{\nu(k)}=\big(-x^*_{0\nu(k)}-x^*_{00},-x^*_{11},\ldots,-x^*_{
\nu(k)-1,\nu(k)-1},0,-a^*_{0\nu(k)}-a^*_{00},-a^*_{11},\ldots,-a^*_{\nu(k)-1\,
\nu(k)-1},0,\\
&-b^*_{0\nu(k)}-b^*_{00},-b^*_{11},\ldots,-b^*_{\nu(k)-1\,\nu(k)-1},0,-u^*_{
0\nu(k)}-u^*_{00},\ldots,-u^*_{\nu(k)-1\,\nu(k)-1},-X^*_{00},\ldots,-X^*_{
\nu(k)-1\,\nu(k)-1},0,\ldots,0\big).
\end{aligned}
\end{eqnarray*}
Next we represent the right-hand side of the inclusion in \eqref{70} by
\begin{eqnarray*}
\lm^k\partial\phi_0(\oz)+\sum_{i=1}^{m}\xi_i^k\nabla
q_i(\oz)+\sum_{j=0}^{\nu(k)-1}\al_{ij}^{1k}\nabla
l^1_{ij}(\oz)+\sum_{j=0}^{\nu(k)-1}\al_{ij}^{2k}\nabla
l^2_{ij}(\oz)+\sum_{j=0}^{\nu(k)-1}\gra g_j(\oz)^*p_{j+1}^k
\end{eqnarray*}
with the complementary slackness conditions
\begin{eqnarray*}
\xi_i^k\(\la a_{i\nu(k)}^k,x_{\nu(k)}^k\ra-b_{i\nu(k)}^k\)=0,\;i=1,\ldots,m.
\end{eqnarray*}
Unifying the above representations and denoting
\begin{equation*}
\rho_j(\oz):=\int_{t_j^k}^{t_{j+1}^k}\n\(\frac{\ox_{j+1}^k-\ox_j^k}{h^k_j}-\dot{
\ox}(t),\frac{\oa_{j+1}^k-\oa_j^k}{h^k_j}-\dot{\oa}(t),\frac{\ob_{j+1}^k-\ob_j^k
}{h^k_j}-\dot{\ob}(t),\ou_j^k(t)-\ou(t)\)
\en^2dt,
\end{equation*}
we arrive at the following relationships:
$$\(\sum_{i=1}^{m}\xi_i^k\nabla
q_i(\oz)\)_{(x_{\nu(k)},a_{\nu(k)},b_{\nu(k)},u_{\nu(k)})}=\(\sum_{i=1}^{m}
\xi_i^k\oa_{ik}^k,\[\xi^k,\rep_m(\ox_{\nu(k)}^k)\],-\xi^k,0\),$$
$$\(\sum_{j=0}^{\nu(k)}\sum^m_{i=1}\al_{ij}^{1k}\nabla
l^1_{ij}(\oz)\)_{(a_j)}=2\[\al_{j}^{1k},\oa_{j}^k\],\;\;j=0,\ldots,\nu(k)-1,$$
$$\(\sum_{j=0}^{\nu(k)}\sum^m_{i=1}\al_{ij}^{2k}\nabla
l^2_{ij}(\oz)\)_{(a_j)}=2\[\al_{j}^{2k},\oa_{j}^k\],\;\;j=0,\ldots,\nu(k)-1,$$
$$\(\sum_{j=0}^{\nu(k)-1}\nabla g_j(\oz)^*p_{j+1}^k\)_{(x_j,a_j,b_j)}=\left\{
\begin{array}{llll}
-p_{1}^{k}&\textrm{ if }& j=0\\[1ex]
p_{j}^{k}-p_{j+1}^{k}&\textrm{ if }& j=1,\ldots,\nu(k)-1\\[1ex]
p_{\nu(k)}^{k}&\textrm{ if }& j=\nu(k)
\end{array}
\right.,
$$
\begin{eqnarray*}
\begin{aligned}
&\(\sum_{j=0}^{\nu(k)-1}\nabla
g_j(\oz)^*p_{j+1}^k\)_{(X_j,A_j,B_j)}=(-h^k_0p_{1}^{xk},-h^k_1p_{2}^{xk},\ldots,
-h^k_{\nu(k)-1}p_{\nu(k)}^{xk},\\
&-h^k_0p_{1}^{ak},-h^k_1p_{2}^{ak},\ldots,-h^k_{\nu(k)-1}p_{\nu(k)}^{ak},
-h^k_0p_{1}^{bk},-h^k_1p_{2}^{bk},\ldots,-h^k_{\nu(k)-1}p_{\nu(k)}^{bk}),
\end{aligned}
\end{eqnarray*}
$$
\partial\phi_0(\oz)\subset\partial\vph(\ox_{\nu(k)}^k)+\sum_{j=0}^{\nu(k)-1}
h^k_j\partial\ell\(\ox_j^k,\oa_j^k,\ob_j^k,\ou_j^k,\oX_j^k,\oA_j^k,
\oB_j^k\)+\frac{1}{2}\sum_{j=0}^{\nu(k)-1}\gra \rho_j(\oz).$$
Furthermore, the set $\lm^k\partial\phi_0(\oz)$ is represented as the collection
of vectors
$$
\begin{aligned}
&\lm^k(h^k_0w_{0}^{xk},h^k_1w_{1}^{xk},\ldots,h^k_{\nu(k)-1}kw_{\nu(k)-1}^{xk},
v_{\nu(k)}^k,h^k_0w_{0}^{ak},h^k_1w_{1}^{ak},\ldots,h^k_{\nu(k)-1}w_{\nu(k)-1}^{
ak},0,\\
&h^k_0w_{0}^{bk},h^k_1w_{1}^{bk},\ldots,h^k_{\nu(k)-1}w_{\nu(k)-1}^{bk},0,\th_{0
}^{uk}+h^k_0w_{0}^{uk},\th_{1}^{uk}+h^k_1w_{1}^{uk},\ldots,\th_{\nu(k)-1}^{uk}
+h^k_{\nu(k)-1}w_{\nu(k)-1}^{uk},\\
&\th_{0}^{Xk}+h^k_0v_{0}^{xk},\th_{1}^{Xk}+h^k_1v_{1}^{xk},\ldots,\th_{\nu(k)-1}
^{Xk}+h^k_{\nu(k)-1}v_{\nu(k)-1}^{xk},\th_{0}^{Ak}+h^k_0v_{0}^{ak},\th_{1}^{Ak}
+h^k_1v_{1}^{ak},\ldots,
\th_{\nu(k)-1}^{Ak}+h^k_{\nu(k)-1}v_{\nu(k)-1}^{ak},\\
&\th_{0}^{Bk}+h^k_0v_{0}^{bk},\th_{1}^{Bk}+h^k_1v_{1}^{bk},\ldots,\th_{\nu(k)-1}
^{Bk}+h^k_{\nu(k)-1}v_{\nu(k)-1}^{bk}),
\end{aligned}
$$
where the components above are such that
$$
\begin{aligned}
&v_{\nu(k)}^k\in\partial\vph(\ox_{\nu(k)}^k),\;\mbox{ and}\\
&\(w_{j}^{xk},w_{j}^{ak},w_{j}^{bk},w_{j}^{uk},v_{j}^{xk},v_{j}^{ak},v_{j}^{bk}
\)\in\partial\ell\(\ox_j^k,\oa_j^k,\ob_j^k,\ou_j^k,\frac{\ox_{j+1}^k-\ox_j^k}{
h^k_j},\frac{\oa_{j+1}^k-\oa_j^k}{h^k_j},
\frac{\ob_{j+1}^k-\ob_j^k}{h^k_j}\),\\
&\left(\th_{j}^{uk},\th_{j}^{Xk},\th_{j}^{Ak},\th_{j}^{Bk}\right):=\\
&
\left(\int_{t_j^k}^{t_{j+1}^k}\(\ou_j^k-\ou(t)\)dt,\int_{t_j^k}^{t_{j+1}^k}
\(\frac{\ox_{j+1}^k-\ox_j^k}{h^k_j}-\dot{\ox}(t)\)dt,\int_{t_j^k}^{t_{j+1}^k}
\(\frac{\oa_{j+1}^k-\oa_j^k}{h^k_j}-\dot{\oa}(t)\)dt,\int_{t_j^k}^{t_{j+1}^k}
\(\frac{\ob_{j+1}^k-\ob_j^k}{h^k_j}-\dot{\ob}(t)\)dt\right)
\end{aligned}
$$
for $j=0,\ldots,\nu(k)-1$. Unifying all of this gives us the conditions
\begin{equation}\label{daux}
-x^*_{00}-x^*_{0\nu(k)}=\lm^kh^k_0w_{0}^{xk}-p_{1}^{xk},
\end{equation}
\begin{equation}\label{e:5.21*x}
-x^*_{jj}=\lm^kh^k_jw_{j}^{xk}+p_{j}^{xk}-p_{j+1}^{xk},\;j=1,\ldots,\nu(k)-1,
\end{equation}
\begin{equation}\label{e:5.24*x}
0=\lm^kv_{\nu(k)}^k+p_{\nu(k)}^{xk}+\sum_{i=1}^{m}\xi_i^k\oa_{ik}^k,\textrm{
where }\;v_{\nu(k)}^k\in\partial\vph (\ox_{\nu(k)}^{k}),
\end{equation}
\begin{equation}\label{daua}
-a^*_{00}-a^*_{0\nu(k)}=\lm^kh^k_0w_{0}^{ak}+2\[\al_{0}^{1k}+\al_{0}^{2k},\oa_{0
}^k\]-p_{1}^{ak},\;\;i=1,\ldots,m,
\end{equation}
\begin{equation}\label{e:5.21*a}
-a^*_{jj}=\lm^kh^k_jw_{j}^{ak}+2\[\al_{j}^{1k}+\al_{j}^{2k},\oa_{j}^k\]+p_{j}^{
ak}-p_{j+1}^{ak},\;\;i=1,\ldots,m,\;\;j=1,\ldots,\nu(k)-1,
\end{equation}
\begin{equation}\label{e:5.24*a}
0=2\(\al_{\nu(k)}^{1k}+\al_{\nu(k)}^{2k}\)\oa_{i\nu(k)}^k+p_{\nu(k)}^{ak}+\[
\xi^k,\rep_m(\ox_{\nu(k)}^k)\],\;\;i=1,\ldots,m,
\end{equation}
\begin{equation}\label{daub}
-b^*_{00}-b^*_{0\nu(k)}=\lm^kh^k_0w_{0}^{bk}-p_{1}^{bk},
\end{equation}
\begin{equation}\label{e:5.21*b}
-b^*_{jj}=\lm^kh^k_jw_{j}^{bk}+p_{j}^{bk}-p_{j+1}^{bk},\;\;j=1,\ldots,\nu(k)-1,
\end{equation}
\begin{equation}\label{e:5.24*b}
0=p_{\nu(k)}^{bk}-\xi^k,
\end{equation}
\begin{equation}\label{dauu}
-u^*_{00}=\lm^k\th_{0}^{uk}+\lm^kh^k_0w_{0}^{uk},
\end{equation}
\begin{equation}\label{e:5.22*}
-u^*_{jj}=\lm^k\th_{j}^{uk}+\lm^kh^k_jw_{j}^{uk},\;\;j=1,\ldots,\nu(k)-1,
\end{equation}
\begin{equation}\label{e:5.23*X}
-X^*_{jj}=\lm^k\th_{j}^{Xk}+\lm^kh_kv_{j}^{xk}-h^k_jp_{j+1}^{xk},\;\;j=0,\ldots,
\nu(k)-1,
\end{equation}
\begin{equation}\label{e:5.23*A}
0=\lm^k\th_{j}^{Ak}+\lm^kh_kv_{j}^{ak}-h^k_jp_{j+1}^{ak},\;\;j=0,\ldots,\nu(k)-1
,
\end{equation}
\begin{equation}\label{e:5.23*B}
0=\lm^k\th_{j}^{Bk}+\lm^kh_kv_{j}^{bk}-h^k_jp_{j+1}^{bk},\;\;j=0,\ldots,
\nu(k)-1.
\end{equation}

Now we are ready to justify all the necessary optimality conditions claimed in
this theorem. First observe that \eqref{mutx}, \eqref{muta}, and \eqref{mutb}
follow from \eqref{e:5.24*x}, \eqref{e:5.24*a}, and \eqref{e:5.24*b},
respectively. Next let us extend each vector $p^k$ by adding the zero component
$p_0^k:=\(x^*_{0\nu(k)},a^*_{0\nu(k)},b^*_{0\nu(k)},u^*_{0\nu(k)}\)$. It follows
from the relationships in \eqref{e:5.21*x}, \eqref{e:5.21*a}, \eqref{e:5.21*b},
\eqref{e:5.23*X}, \eqref{e:5.23*A}, and \eqref{e:5.23*B} that
$$
\begin{aligned}
\frac{x^*_{jj}}{h^k_j}&=\frac{p_{j+1}^{xk}-p_{j}^{xk}}{h^k_j}-\lm^kw_{j}^{xk},\\
\frac{a^*_{jj}}{h^k_j}&=\frac{p_{j+1}^{ak}-p_{j}^{ak}}{h^k_j}-\lm^kw_{j}^{ak}
-\frac{2}{h^k_j}\(\al_{j}^{1k}+\al_{j}^{2k}\)\oa_{ij}^k,\\
\frac{b^*_{jj}}{h^k_j}&=\frac{p_{j+1}^{bk}-p_{j}^{bk}}{h^k_j}-\lm^kw_{j}^{bk},\\
\frac{u^*_{jj}}{h^k_j}&=-\frac{1}{h^k_j}\lm^k\th_{j}^{uk}-\lm^kw_{j}^{uk},\\
\frac{X^*_{jj}}{h^k_j}&=-\frac{1}{h^k_j}\lm^k\th_{j}^{Xk}+p_{j+1}^{xk}-\lm^kv_{j
}^{xk},\\
0&=-\frac{1}{h^k_j}\lm^k\th_{j}^{Ak}+p_{j+1}^{ak}-\lm^kv_{j}^{ak},\\
0&=-\frac{1}{h^k_j}\lm^k\th_{j}^{Bk}+p_{j+1}^{bk}-\lm^kv_{j}^{bk}.
\end{aligned}
$$
Substituting this into the left-hand side of \eqref{e:5.18*} and taking into
account the equalities in \eqref{71+}--\eqref{71l2}, \eqref{e:5.24*x},
\eqref{e:5.24*a}, and \eqref{e:5.24*b} justify the claims made in
\eqref{xi}--\eqref{e:5.10*}.

To verify finally the nontriviality condition \eqref{e:5.8*}, suppose on the
contrary that
$\lm^k=0,\,\xi^k=0,\,\al^{1k}+\al^{2k}=0,\,p_{j}^{xk}=0,\,p_{j}^{ak}=0,\,p_{j}^{
bk}=0,\,\psi^k=0$ for all $j=0,\ldots,\nu(k)-1$, which yields in turn
$x^*_{0k}=p_{0}^{xk}=0,\,a^*_{0k}=p_{0}^{ak}=0$, and $b^*_{0k}=p_{0}^{bk}=0$.
Then it follows from \eqref{e:5.24*x}, \eqref{e:5.24*a}, and \eqref{e:5.24*b}
that $\(p_{\nu(k)}^{xk},p_{\nu(k)}^{ak},p_{\nu(k)}^{bk}\)=0$, and hence
$\(p_{j}^{xk},p_{j}^{ak},p_{j}^{bk}\)=0$, for all $j=0,\ldots,\nu(k)$. We see
also that the conditions in \eqref{daux}, \eqref{e:5.21*x}, \eqref{daua},
\eqref{e:5.21*a}, \eqref{daub}, \eqref{e:5.21*b}, \eqref{dauu}, and
\eqref{e:5.22*} imply that $\(x^*_{jj},a^*_{jj},b^*_{jj},u^*_{ij}\)=0$ for all
$j=0,\ldots,\nu(k)-1$. In addition, it follows from \eqref{e:5.23*X},
\eqref{e:5.23*A}, and \eqref{e:5.23*B} that
$X^*_{jj}=0,\,A^*_{jj}=0,\,B^*_{jj}=0$ for all $j=0,\ldots,\nu(k)-1$.
Furthermore, all the components of $z^*_j$ different from $(x^*_{jj},
a^*_{jj},b^*_{jj},u^*_{jj},X^*_{jj},A^*_{jj},B^*_{jj})$ are clearly zero for
$j=0,\ldots,\nu(k)-1$, and hence $z^*_{j}=0$ for $j=0,\ldots,\nu(k)-1$. We
similarly conclude that $z^*_{k}=0$, since  $x^*_{0k}=p_{0}^{xk}=0$ while all
the other components of this vector obviously reduce to zero. Thus $z^*_j=0$ for
all $j=0,\ldots,\nu(k)$, which violates the nontriviality condition for $(MP)$
and completes the proof of the theorem. $\h$

Our next theorem provides verifiable necessary optimality conditions for
solutions $(\ox^k,\oa^k,\ob^k,\ou^k)$ to problems $(P_k)$ that strongly
approximate the given r.i.l.m.\ $(\ox,\oa,\ob,\ou)$ for the original sweeping
control problem $(P)$. The proof is based on the results of Theorem~\ref{Th8}
and the second-order calculations from Theorem~\ref{Th7}.\vspace*{-0.1in}

\begin{theorem}{\bf(optimality conditions for discretized sweeping processes via
their initial data)}\label{Th9} Let $(\ox^k,\oa^k,\ob^k,\ou^k)$ be an optimal
solution to problem $(P_k)$ under the notation and assumptions of
Theorem~{\rm\ref{Th8}} for each fixed index $k\in\N$. Then there exist dual
elements $(\lm^k,\al^{1k},\al^{2k},\psi^k,p^k)$ as in Theorem~{\rm\ref{Th8}}
together with vectors $\eta_j^k\in\R^m_+$ as $j=0,\ldots,\nu(k)-1$ and
$\gg_j^k\in\R^m$ as $j=0,\ldots,\nu(k)-1$ satisfying the following
conditions:\\[2ex]
$\bullet$ The {\sc primal arc representation:}
\begin{equation}\label{87}
-\frac{\ox_{j+1}^k-\ox_j^k}{h^k_j}+g(\ox_j^k,\ou_j^k)=\sum_{i=1}^{m}\eta_{ij}
^k\oa_{ij}^k,\quad j=0,\ldots,\nu(k)-1.
\end{equation}
$\bullet$ The {\sc adjoint dynamic relationships:}
\begin{equation}\label{conx}
\frac{p_{j+1}^{xk}-p_{j}^{xk}}{h^k_j}-\lm^kw_{j}^{xk}\in\nabla_x
g(\ox_j^k,\ou_j^k)^*\Big(\frac{1}{h^k_j}\lm^k\th_{j}^{Xk}+\lm^k
v_{j}^{xk}-p_{j+1}^{xk}\Big)+\sum_{i=1}^{m}\gg_{ij}^k\oa_{ij}^k,
\end{equation}
\begin{equation}\label{cona}
\dfrac{p_{j+1}^{ak}-p_{j}^{ak}}{h^k_j}-\lm^kw_{j}^{ak}-\dfrac{2}{h^k_j}\[\al_{j}
^{1k}+\al_{j}^{2k},\oa_{j}^k\]=\[\gg_{j}^k,\rep_m(\ox_j^k)\]+\Big[\eta_j^k,
\rep_m\Big(-\dfrac{1}{h^k_j}\lm^k\th_{j}^{Xk}-\lm^k
v_{j}^{xk}+p_{j+1}^{xk}\Big)\Big],
\end{equation}
\begin{equation}\label{conb}
\frac{p_{j+1}^{bk}-p_{j}^{bk}}{h^k_j}-\lm^kw_{j}^{bk}=-\gg_{j}^k,\quad
j=0,\ldots,\nu(k)-1,
\end{equation}
where the components of the vectors $\gg^k_j$ are such that
\begin{equation*}
\begin{cases}
\gg^k_{ij}=0\;\;\mbox{if}\;\;\la\oa^k_{ij},\ox^k_j\ra<\ob^k_{ij},\;\;\mbox{or}
\;\;\eta^k_{ij}=0\;\mbox{ and
}\;\la\oa_{ij}^k,-\dfrac{1}{h^k_j}\lm^k\th_{j}^{Xk}-\lm^k
v_{j}^{xk}+p_{j+1}^{xk}\ra<0,\\
\gg^k_{ij}\ge
0\;\;\mbox{if}\;\;\la\oa^k_{ij},\ox^k_j\ra=\ob^k_{ij},\;\eta^k_{ij}=0,\;\mbox{
and }\;\la\oa_{ij}^k,-\dfrac{1}{h^k_j}\lm^k\th_{j}^{Xk}-\lm^k
v_{j}^{xk}+p_{j+1}^{xk}\ra>0,\\
\gg^k_{ij}\in\R\;\;\mbox{if}\;\;\eta^k_{ij}>0\;\mbox{ and
}\;\la\oa_{ij}^k,-\dfrac{1}{h^k_j}\lm^k\th_{j}^{Xk}-\lm^k
v_{j}^{xk}+p_{j+1}^{xk}\ra=0
\end{cases}
\end{equation*}
for the indices $j=0,\ldots,\nu(k)-1$ and $i=1,\ldots,m$.\\[2ex]
$\bullet$ The {\sc local maximum principle:}
\begin{equation}\label{cony}
\psi^k_j\in N(\ou^k_j;U)\;\mbox{ with
}\;-\dfrac{1}{h^k_j}\psi_j^k-\dfrac{1}{h^k_j}\lm^k\th_{j}^{uk}-\lm^kw_{j}^{uk}
\in\nabla_u g(\ox_j^k,\ou_j^k)^*\Big(\frac{1}{h^k_j}\lm^k\th_{j}^{Xk}+\lm^k
v_{j}^{xk}-p_{j+1}^{xk}\Big)
\end{equation}
for $j=0,\ldots,\nu(k)-1$, where the subgradients
$(w_{j}^{xk},w_{j}^{ak},w_{j}^{bk},w_{j}^{uk},v_{j}^{xk},v_{j}^{ak},v_{j}^{bk})$
are taken from \eqref{subcol}. If furthermore the normal cone $N(\ou^k_j;U)$ is
tangentially generated, i.e.,
\begin{equation*}
N(\ou^k_j;U)=T^*(\ou^k_j;U):=\big\{v\in\R^d\big|\;\la v,u\ra\le 0\;\mbox{ for
all }\;u\in T(\ou^k_j;U)\big\},
\end{equation*}
for some tangent cone $T(\ou^k_j;U)$, then the first inclusion in \eqref{cony}
is written as
\begin{equation}\label{lmp}
\la\psi^k_j,\ou^k_j\ra=\max_{u\in T(\ou^k_j;U)}\la\psi^k_j,u\ra,\quad
j=0,\ldots,\nu(k)-1,
\end{equation}
which reduces to the {\sc global maximum principle}
\begin{equation}\label{gmp}
\langle\psi^k_j,\ou^k_j\rangle=\max_{u\in U}\la\psi^k_j,u\ra,\quad
j=0,\ldots,\nu(k)-1,
\end{equation}
provided that the control set $U$ is convex.\\[2ex]
$\bullet$ The {\sc transversality conditions} at the right endpoint:
\begin{equation}\label{nmutx}
-p_{\nu(k)}^{xk}\in\lm^k\partial\vph(\ox_{\nu(k)}^k)+\sum_{i=1}^{m}\eta_{i\nu(k)
}^k\oa_{i\nu(k)}^k,
\end{equation}
\begin{equation}\label{nmuta}
p_{\nu(k)}^{ak}=-2\[\al_{\nu(k)}^{1k}+\al_{\nu(k)}^{2k},\oa_{i\nu(k)}^k\]-\[
\eta_{\nu(k)}^k,\rep_m(\ox_{\nu(k)}^k)\],
\end{equation}
\begin{equation}\label{nmutb}
p_{i\nu(k)}^{bk}=\eta^k_{i\nu(k)}\ge 0,\;\la \oa^{k}_{i\nu(k)},\ox^k_{\nu(k)}\ra
< \ob^k_{i\nu(k)}\sr p_{i\nu(k)}^{bk}=0\;\;\textrm{for}\;\;i=1,\ldots,m
\end{equation}
with dual vectors $\al_{\nu(k)}^{1k}$ and $\al_{\nu(k)}^{2k}$ satisfying
\begin{equation}\label{t:7.20}
\al_{i\nu(k)}^{1k}\in N_{\[0,1+\dd_k\]}(\|\oa_{i\nu(k)}^k\|)\;\mbox{ and
}\;\al_{i\nu(k)}^{2k}\in
N_{\[1-\dd_k,\infty\]}(\|\oa_{i\nu(k)}^k\|),\;i=1,\ldots,m,
\end{equation}
where the normal cone to the convex sets is explicitly expressed in form
\eqref{NC}.\\[2ex]
$\bullet$ The {\sc complementarity slackness conditions:}
\begin{equation} \label{eta}
\[\la a_{ij}^k,\ox_j^k\ra < \ob_{ij}^k\]\sr\eta_{ij}^k=0,
\end{equation}
\begin{equation}\label{eta1}
\[\la\oa_{i\nu(k)}^k,\ox_{\nu(k)}^k\ra<\ob_{i\nu(k)}^k\]\sr\eta_{i\nu(k)}^k=0,
\end{equation}
\begin{equation}\label{96}
\eta_{ij}^k>0\sr\[\la\oa_{ij}^k,-\frac{1}{h^k_j}\lm^k\th_{j}^{Xk}-\lm^kv_{j}^{xk
}+p_{j+1}^{xk}\ra=0\]
\end{equation}
for all the indices $j=0,\ldots,\nu(k)-1$ and $i=1,\ldots,m$.\\[2ex]
$\bullet$ The {\sc nontriviality conditions:}
\begin{equation}\label{ntc}
\lm^k+\n\al^{1k}+\al^{2k}\en+\n\eta_{\nu(k)}^k\en +\sum_{j=0}^{\nu(k)-1}\n
p_{j}^{xk}\en +\n p_{0}^{ak}\en+\n p_{0}^{bk}\en+ \n\psi^k\en\ne 0,
\end{equation}
\eq
\label{ntc1}
\lm^k+\n\al^{1k}+\al^{2k}\en+\n\gg^k\en\ne 0.
\eeq
\end{theorem}\vspace*{-0.1in}
{\bf Proof}. It follows from condition \eqref{e:5.10*} of Theorem~\ref{Th8} and
the coderivative definition \eqref{coderivative} that
\begin{equation*}
\begin{aligned}
&\bigg(\dfrac{p_{j+1}^{xk}-p_{j}^{xk}}{h^k_j}-\lm^kw_{j}^{xk},\dfrac{p_{j+1}^{ak
}-p_{j}^{ak}}{h^k_j}-\lm^kw_{j}^{ak}-\dfrac{2}{h^k_j}\(\al_{j}^{1k}+\al_{j}^{2k}
\)\oa_{ij}^k,\dfrac{p_{j+1}^{bk}-p_{j}^{bk}}
{h^k_j}-\lm^kw_{j}^{bk},\\&-\dfrac{1}{h^k_j}\lm^k\th_{j}^{uk}-\lm^kw_{j}^{uk}
-\dfrac{1}{h^k_j}\psi_j^k\bigg)\in
D^*F\(\ox_j^k,\oa_j^k,\ob_j^k,\ou_j^k,-\frac{\ox_{j+1}^k-\ox_j^k}{h^k_j}
\)\(-\frac{1}{h^k_j}\lm^k\th_{j}^{Xk}-\lm^k v_{j}^{xk}+p_{j+1}^{xk}\)
\end{aligned}
\end{equation*}
for all $j=0,\ldots,\nu(k)-1,\;i=1,\ldots,m$. Using the inclusion
\begin{equation*}
-\frac{\ox_{j+1}^k-\ox_j^k}{h^k_j}+g(\ox_j^k,\ou_j^k)\in
G(\ox_j^k,\oa_j^k,\ob_j^k)
\end{equation*}
via the normal cone mapping $G$ from \eqref{G} and employing the PLICQ property
of the vectors $\nn \oa^k_i|\;i\in I(\ox^k,\oa^k,\ob^k)\hnn$ give us a unique
vector $\eta_j^k\in\R^m_+$ such that for all $i=1,\ldots,m$ we have
\begin{equation*}\label{h:5.39}
\sum_{i=1}^m\eta_{ij}^k\oa_{ij}^k=-\frac{\ox_{j+1}^k-\ox_j^k}{h^k_j}+g(\ox_j^k,
\ou_j^k)\;\textrm{ with }\;\eta_{ij}^k\in
N_{\R_-}\Big(\la\oa_{ij}^k,\ox_j^k\ra-\ob_{ij}^k\Big),\;j=0,\ldots,\nu(k)-1,
\end{equation*}
which verifies the implications in \eqref{87} and \eqref{eta}. Applying now the
coderivative upper estimate \eqref{c57} from Theorem~\ref{Th7} with
$x:=\ox_j^k$, $a:=\oa_j^k$, $b:=\ob_j^k$, $u:=\ou_j^k$,
$w:=-\dfrac{\ox_{j+1}^k-\ox_j^k}{h^k_j}$, and
$y:=-\dfrac{1}{h_k}\lm^k\th_{j}^{Xk}-\lm^k v_{j}^{xk}+p_{j+1}^{xk}$ as
$j=0,\ldots,\nu(k)-1$ shows that $\gg_j^k\in\R^m$ and that the relationships
$$
\(\begin{matrix}
\dfrac{p_{j+1}^{xk}-p_{j}^{xk}}{h^k_j}-\lm^kw_{j}^{xk},\dfrac{p_{j+1}^{ak}-p_{j}
^{ak}}{h^k_j}-\lm^kw_{j}^{ak}-\dfrac{2}{h^k_j}\[\al_{j}^{1k}+\al_{j}^{2k},\oa_{j
}^k\],\dfrac{p_{j+1}^{bk}-p_{j}^{bk}}{h^k_j}
-\lm^kw_{j}^{bk},\\
-\dfrac{1}{h^k_j}\lm^k\th_{j}^{uk}-\lm^kw_{j}^{uk}-\dfrac{1}{h^k_j}\psi_j^k
\end{matrix}\)
$$
$$
\in\(\begin{matrix}
\disp-\nabla g_x(\ox_j^k,\ou_j^k)^*\(-\frac{1}{h^k_j}\lm^k\th_{j}^{Xk}-\lm^k
v_{j}^{xk}+p_{j+1}^{xk}\)+\sum_{i=1}^{m}\gg_{ij}^k\oa_{ij}^k,\\
\[\gg_{j}^k,\rep_m(\ox_j^k)\]+\[\eta_j^k,\rep_m\(-\dfrac{1}{h^k_j}\lm^k\th_{j}^{
Xk}-\lm^k v_{j}^{xk}+p_{j+1}^{xk}\)\],\\
-\gg_{j}^k,\;\disp-\nabla
g_u(\ox_j^k,\ou_j^k)^*\(-\frac{1}{h^k_j}\lm^k\th_{j}^{Xk}-\lm^k
v_{j}^{xk}+p_{j+1}^{xk}\)
\end{matrix}\)
$$
are satisfied, where $\psi_{j}^{k}\in N(\ou^k_j;U)$ for all
$j=0,\ldots,\nu(k)-1$, and where the components $\gg^k_{ij}$ of the vectors
$\gg_j^k\in\R^m$ as $i=1,\ldots,m$ are taken from
\begin{equation}\label{congg}
\gg^k_{ij}\in
D^*N_{\R_-}\(\la\oa_{ij}^k,\ox_j^k\ra-\ob_{ij}^k,\eta^k_{ij}\)\(\la\oa_{ij}^k,
-\dfrac{1}{h^k_j}\lm^k\th_{j}^{Xk}-\lm^k v_{j}^{xk}+p_{j+1}^{xk}\ra\).
\end{equation}
The obtained relationships together with the direct calculation of the
coderivative $D^*N_{\R_-}$ in \eqref{congg} ensure the validity of all the
conditions in \eqref{conx} as well as the inclusion in \eqref{cony}. The latter
together with \eqref{psi} constitutes an appropriate version of the (linearized)
local maximum principle for nonconvex discrete-time systems. It immediately
gives us the local maximality condition \eqref{lmp} in the case of tangentially
generated normals, which surely holds for the class of normally regular sets
$U$; see, e.g., \cite{m-book,rw}. The global form of the discrete maximum
principle in \eqref{gmp} is a direct consequence of \eqref{cony} and the normal
cone representation \eqref{NC} for convex sets. Furthermore, conditions
\eqref{nmutx}, \eqref{nmuta}, and \eqref{nmutb} clearly follow from
\eqref{mutx}, \eqref{muta}, and \eqref{mutb} due to \eqref{xi}.

Defining now $\eta_{\nu(k)}^k:=\xi^k$ via $\xi^k$ from the statement of
Theorem~\ref{Th8} yields $\eta_j^k\in\R^m_+$ for $j=0,\ldots,\nu(k)$ and allows
us to deduce the nontriviality condition \eqref{ntc} from that in \eqref{e:5.8*}
and also the transversality conditions in \eqref{nmutx}--\eqref{nmutb} from
those in \eqref{mutx}-- \eqref{mutb}. Implication \eqref{eta1} is a direct
consequence of \eqref{xi} and the definition of $\eta_{\nu(k)}^k$. Observing
that \eqref{96} follows from the fact that
\begin{equation*}
-\dfrac{1}{h^k_j}\lm^k\th_{j}^{Xk}-\lm^k
v_{j}^{xk}+p_{j+1}^{xk}\in\disp\bigcap_{\nn
i|\;\eta^k_{ij}>0\hnn}(\oa^k_{ij})^\perp,
\end{equation*}
we get from \eqref{al1} and \eqref{al2} that both inclusions in \eqref{t:7.20}
hold.

It remains to verify the nontriviality condition \eqref{ntc1}. Suppose on the
contrary that $\lm^k=0,\; \al^{1k}+\al^{2k}=0$, and $\gg^k=0$. We deduce from
\eqref{e:dac5} that $p^{ak}_{\nu(k)}=0$ and $p^{bk}_{\nu(k)}=0$, which clearly
yield $\eta^k_{\nu(k)}=p^{bk}_{\nu(k)}=0$. Then it follows from \eqref{nmutx}
that $p^{xk}_{\nu(k)}=0$, and thus $\(p^{xk}_j,p^{ak}_j\)=(0,0)$ for all
$j=0,\ldots,\nu(k)-1$ by \eqref{conx} and \eqref{cona}. This implies that
$\psi^k=0$ by \eqref{cony}. Using finally \eqref{conb} tells us that
$p^{bk}_0=0$. It means that \eqref{ntc} is violated, which is a contradiction
that justifies the validity of \eqref{ntc1} and therefore completes the proof of
the theorem. $\h$\vspace*{-0.2in}

\section{Numerical Illustration}\label{exam}
\setcounter{equation}{0}\vspace*{-0.1in}

In this section we present a nontrivial example illustrating the application of
the obtained results to solve the sweeping optimal control problem $(P)$. We
consider this problem
with the following data, where the $a$-components and $b$-components of controls
are fixed, and only the $u$-components are used for optimization:
\begin{equation}\label{data}
\left\{
\begin{array}{ll}
n=2,\;m=1,\;T=1,\;x_0=\(\frac{3}{2},1\),\;a=\(-\frac{1}{\sqrt5},-\frac{2}{\sqrt5
}\),\;b=-\frac{2}{\sqrt5},\\[1ex]
g(x,u):=u,\;\;\vph(x):=x_1+x_2,\;\;\ell(t,x,a,b,u,\dot x,\dot a,\dot
b):=\frac{1}{2}u^2_1+u^2_2,\\[1ex]
U:=[-1,1]\times[-1,1].
\end{array}\right.
\end{equation}
The set $C(t)$ in the sweeping inclusion \eqref{SP} is described now by
$$
C(t):=\nn(x_1,x_2)\in\R^2\big|\;x_1+2x_2\ge 2\hnn\;\mbox{ for all }\;t\in[0,1].
$$
In what follows we are going to show that applying the optimality conditions of
Theorem~\ref{Th9} allows us to find optimal solutions to problems $(P_k)$, which
can be viewed as (sub)optimal solutions to the original sweeping control problem
$(P)$. For simplicity and convenience, we consider only the case where $k=2$
(and drop below this superscript), while the calculations are similar for any
natural number $k$.

\begin{figure}[ht]
\centering
\includegraphics[scale=0.4]{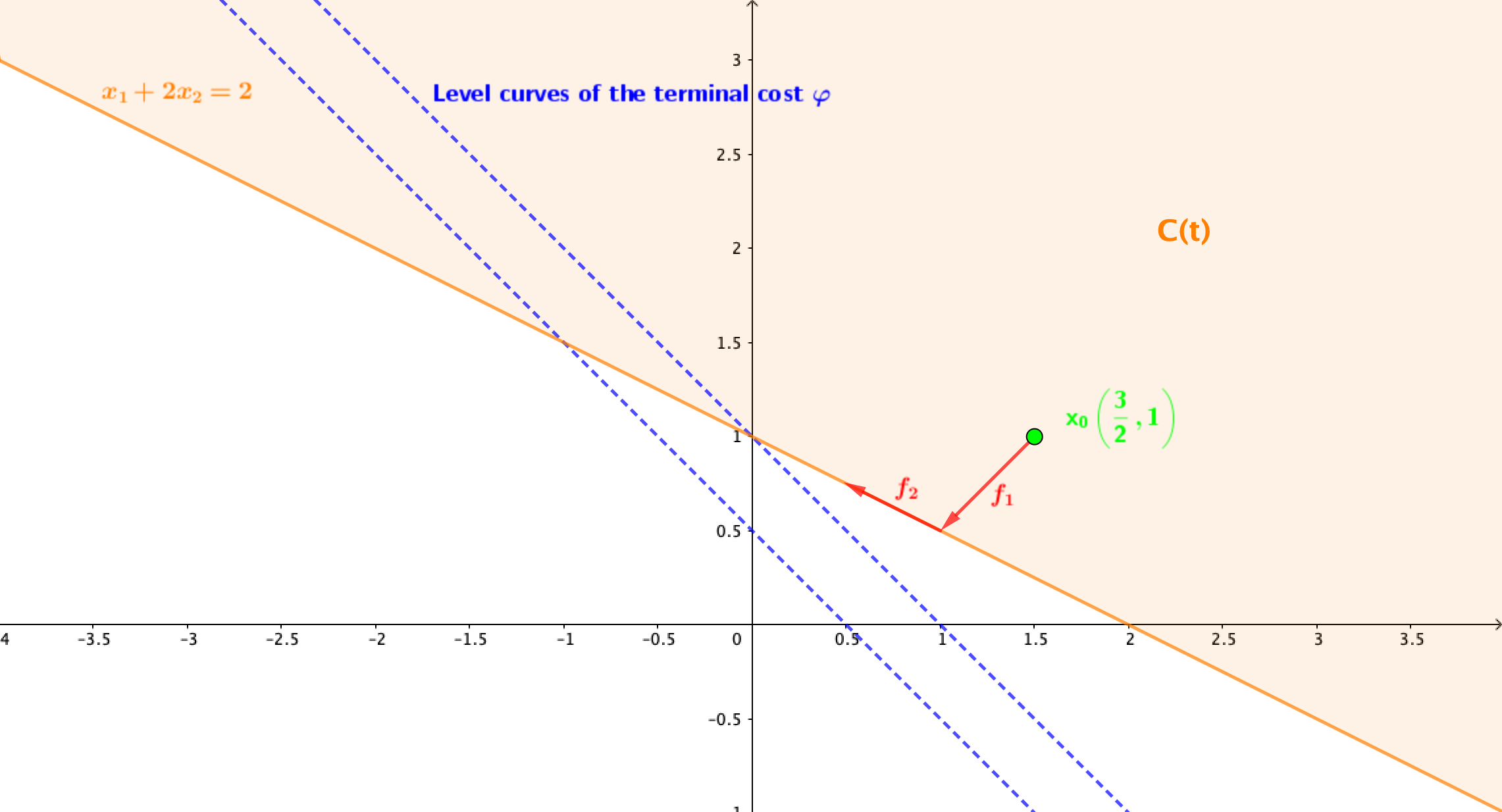}
\caption{Dynamics of the controlled sweeping process.}
\label{Fig1}
\end{figure}

It is easy to see that all the assumptions of Theorem~\ref{Th9} are satisfied
for \eqref{data}. Employing the obtained necessary optimality conditions in this
setting gives us dual elements $\lm\ge 0,\;\eta_j\ge 0,\;\gg_j\in\R$,
$\al^1_j,\al^2_j\in\R$, $\psi_j\in\R^2$, and $(p^x_j,p^a_j,p^b_j)\in\R^5$,
$(w^x_j,w^a_j,w^b_j,w^u_j)\in\R^7$, and $(x^x_j,v^a_j,v^b_j)\in\R^5$ as $j=0,1$
satisfying the following relationships, where
$\(\th^u_j,\th^X_j,\th^A_j,\th^B_j\)\approx 0$ as $j=0,1$ due to the established
convergence of optimal solutions:
\begin{enumerate}
\item $\(w^x_j,w^a_j,w^b_j,w^u_j\)=(0,0,0,0,0,\ou_1,2\ou_2)$ for $j=0,1$.
\item $\(v^x_j,v^a_j,v^b_j\)=(0,0,0,0,0,0)$ for $j=0,1$.
\item
$
\dot\ox(t)=\nn
\begin{array}{lll}
\ou(t)+\eta_0(1,2)&\mbox{if}&t\in (0,\frac{1}{2})\\[1ex]
\ou(t)+\eta_1(1,2)&\mbox{if}&t\in(\frac{1}{2},1)
\end{array}
\right.,
$
$\;\;$where$\;$
$
\ou(t)=\nn
\begin{array}{lll}
\ou_0&\mbox{if}&t\in [0,\frac{1}{2})\\[1ex]
\ou_1&\mbox{if}&t\in(\frac{1}{2},1]
\end{array}
\right.
$.
\item
$
\begin{cases}
2\(p^x_{j+1}-p^x_j\)=\gg_j(1,2),\\
2\(p^a_{j+1}-p^a_j\)-4\(-\(\al^1_{j}+\al^2_{j}\)\frac{1}{\sqrt5},-\(\al^1_{j}
+\al^2_{j}\)\frac{2}{\sqrt5}\)=
\(\gg_{j}\ox_{1j},\gg_{j}\ox_{2j}\)+\(\eta_{j}p^x_{1,j+1},\eta_{j}p^x_{2,j+1}\),
\\
2\(p^b_{j+1}-p^b_j \)=-\gg_j\;\mbox{ for }\;j=0,1.
\end{cases}
$\\[1ex]
\item $2\psi_j+\lm\(\ou_{j1},2\ou_{j2}\)=p^x_{j+1}$ for $j=0,1$.
\item $\psi_j\in N\(\ou_j;[-1,1]\times[-1,1]\)$ for $j=0,1$, which is equivalent
to\\[1ex]
$\psi_{1j}\ou_{1j}+\psi_{2j}\ou_{2j}=
\disp\max_{(u_1,u_2)\in[-1,1]\times[-1,1]}\nn
\psi_{1j}u_{1}+\psi_{2j}u_{2}\hnn$.
\item $\ox_{j1}+2\ox_{j2}>2\sr\gg_j=0\;\mbox{ and }\;\eta_j=0$ for $j=0,1$.
\item $\eta_j>0\sr\la(-1,-2),p^x_{j+1}\ra=0$ for $j=0,1$.
\item $\ox_{21}+2\ox_{22}>2\sr\eta_2=0$.
\item
$
\begin{cases}
-p^x_2=\lm(1,1)+\eta_1(1,2),\\
p^a_2=-2\(-\(\al^1_{2}+\al^2_{2}\)\frac{1}{\sqrt5},-\(\al^1_{2}+\al^2_{2}\)\frac
{2}{\sqrt5}\)-\(\eta_{1}\ox_{12},\eta_{1}\ox_{22}\),\\
p^b_2=\eta_1\ge 0.
\end{cases}
$
\item $\al^1_2\in N_{[0,1+\dd_k]}(1),\; \al^2_2\in N_{[1-\dd_k,\infty)}(1)$,
which implies that $\(\al^1_2,\al^2_2\)=(0,0)$.
\item $\lm+\n\al^1+\al^2\en+\n\gg\en>0$.
\end{enumerate}
It then follows from (2) that
$$
\ox(t)=\nn
\begin{array}{ll}
\(\frac{3}{2}+t\ou_{01},1+t\ou_{02}\)&\mbox{ if }\;t\in[0,\frac{1}{2})\\[1ex]
\(\frac{3}{2}+\frac{1}{2}\ou_{01}+\(t-\frac{1}{2}\)\(\ou_{11}+\eta_1\),1+\frac{1
}{2}\ou_{02}+\(t-\frac{1}{2}\)\(\ou_{12}+2\eta_1\))\) &\mbox{  if  }
t\in[\frac{1}{2},1]
\end{array}
\right.
$$
as $\eta_0=0$ due to (6). Let $t^*$ be the time when the moving particle hits
the boundary, i.e., $\ox_1(t^*)+2\ox_2(t^*)=2$. Consequently, we have that
$$
\frac{7}{2}+t^*(\ou_{01}+2\ou_{02})=2\;\mbox{ if }\;t^*<\frac{1}{2},
$$
$$
\frac{7}{2}+\frac{1}{2}(\ou_{01}+2\ou_{02}) +
\(t^*-\frac{1}{2}\)\(\ou_{11}+2\ou_{12}+5\eta_1\)=2\;\mbox{ if }
\;t^*\ge\frac{1}{2}.
$$
When $\ox(\cdot)$ hits the boundary of the set $C$, it would stay there while
pointing in the direction shown in Figure~\ref{Fig1}.
Thus $t^*=\frac{1}{2}$ and so $\ou_{01}+2\ou_{02}=-3$, which implies that
$\ou_{01}=\ou_{02}=-1$ due to $-1\le\ou_{01},\ou_{02}\le 1$. Since the particle
stays on the boundary after $t=t^*=\frac{1}{2}$, we have $\ox_1(t)+2\ox_2(t)=2$
for all $t\in\big(\frac{1}{2},1]$. Therefore
$$
\frac{7}{2}+\frac{1}{2}(\ou_{01}+2\ou_{02})+\(t-\frac{1}{2}\)\(\ou_{11}+2\ou_{12
}+5\eta_1\)=2\;\mbox{ for all }\;t>\frac{1}{2},
$$
which ensures in turn that
\begin{equation}
\label{u-eta}
\ou_{11}+2\ou_{12}+5\eta_1=0.
\end{equation}
The cost functional is thus calculated by
$$
J[\ox,\ou]=\dfrac{\ou^2_{11}}{4}+\dfrac{\ou^2_{12}}{2}+\frac{1}{2}\(\ou_{11}
+\ou_{12}+3\eta_1\)+\frac{9}{4}.
$$

To proceed further, consider the following two cases:\\[1ex]
{\bf Case~1:} $\eta_1=0$, i.e., the normal vector $\eta_1(-1,-2)$ taken from the
normal cone $N(x(t);C(t))$ is not active for $t>\frac{1}{2}$. It then follows
from \eqref{u-eta} that
$\ou_{11}=-2\ou_{12}$. The cost functional in this case is
$$
J[\ox,\ou]=\frac{3}{2}\ou^2_{12}-\frac{1}{2}\ou_{12}+\frac{9}{4},
$$
and it attains the minimum value at $\ou_{12}=\frac{1}{6}$. Thus
$\ou_{11}=-\frac{1}{3}$, and the minimum value of $J$ is $\frac{53}{24}\approx
2.208$.\\[1ex]
{\bf Case~2:} $\eta_1>0$. Using (8) gives us $p^x_{21}+2p^x_{22}=0$. On the
other hand, we obtain from (4) and (5) that
\begin{equation}
\label{u-re}
2\psi_{11}+4\psi_{12}+\lm\(\ou_{11}+4\ou_{12}\)=0.
\end{equation}
Examine now the two possibilities:
\begin{itemize}
\item If either $\psi_{11}\not=0$ or $\psi_{12}\not=0$, it follows from (6) that
either $\ou_{11}=\pm1$ or $\ou_{12}=\pm1$. The possible minimum values of the
cost functional in this case are
collected the following table, where the symbol `$X$' indicates that the control
in question is not admissible:
\begin{center}
\begin{tabular}{|c|c|c|}
\hline
The value of the controls & The minimum value of the cost functional & The value
of $\eta_1$\\[1ex]
\hline
$\ou_{11}=1, \ou_{12}=\frac{1}{10}$  &X  & negative (not admissible)\\[1ex]
\hline
$\ou_{11}=-1, \ou_{12}=\frac{1}{5}$  &3.295  & $\frac{4}{25}$ \\[1ex]
\hline
$\ou_{12}=1, \ou_{11}=-\frac{2}{5}$  &X  & negative (not admissible)\\[1ex]
\hline
$\ou_{12}=-1, \ou_{11}=-\frac{2}{5}$ &3.81 & $\frac{12}{25}$ \\[1ex]
\hline
\end{tabular}
\end{center}

\item If $\psi_1=0$, then $\ou_{11}+4\ou_{12}=0$ by \eqref{u-re} assuming that
$\lm>0$; otherwise we do not have enough information to
proceed. In this case the cost functional is
$$
J[\ox,\ou]=\frac{9}{2}\ou^2_{12}-\frac{9}{5}\ou_{12}+\frac{9}{4},
$$
which achieves the minimum value at $\ou_{12}=\frac{1}{10}$. The minimum cost
value in this case is $\frac{441}{200}=2.205$, the other component of the
control is $\ou_{11}=-\frac{2}{5}$, and the associated value of $\eta_1$ is
$\frac{1}{25}$.
\end{itemize}
The corresponding optimal trajectory of the sweeping process is calculated by
$$
\ox(t)=\nn
\begin{array}{ll}
\(\frac{3}{2}-t,1-t\)&\mbox{ if }\;t\in[0,\frac{1}{2})\\[1ex]
\Big(1-\frac{9}{25}\(t-\frac{1}{2}\),\frac{1}{2}+\frac{9}{50}\(t-\frac{1}{2}
\)\Big)&\mbox{ if }\;t\in[\frac{1}{2},1]
\end{array}
\right..
$$

Note also that we can reduce the cost functional to the function of two
variables $\ou_{11}$ and $\ou_{12}$ by solving \eqref{u-eta} for $\eta_1$ and
arriving then to the following expression:
$$
J[\ox,\ou]=\frac{1}{4}\ou^2_{11}+\frac{1}{2}\ou^2_{12}+\frac{1}{5}\ou_{11}-\frac
{1}{10}\ou_{12}+\frac{9}{4}.
$$
Thus we can treat our problem as minimizing the latter objective function under
the inequality constraints $\eta_1=-\frac{1}{5}\ou_{11}-\frac{2}{5}\ou_{12}\ge0$
and $(\ou_{11},\ou_{12})\in[-1,1]\times[-1,1]$. The solution obtained in this
way agrees with the one calculated above.

Finally, let us summarize in the next remark the mechanism of determining
optimal controls developed above, which is illustrated by
Figure~1.\vspace*{-0.1in}

\begin{remark}\label{exp}
In order to reduce the terminal cost effectively for the controlled sweeping
process in this example, the trajectory $\ox(\cdot)$ should follow the direction
of the vector $(-1,-1)$ (a negative gradient vector of the terminal cost) while
keeping as least energy as possible. When the particle hits the boundary of the
set $C(\cdot)$, it would stay there and point in the direction of the vector
$(-2,1)$ to keep reducing the value of $\ox_1(t)+\ox_2(t)$ till the end of the
process. To make this happen, it is natural (at the first glance) to push the
particle horizontally as twice as vertically so that it points in the direction
of $(-2,1)$. Roughly speaking, we tend to use more energy for the first
component $\ou_{11}$ of the control $\ou(t)$ than for the other component in
order to force the particle to point in the desired direction. That is due to
the form of the running cost as in \eqref{data}. If the normal vector
$\eta_1(-1,-2)$ generated from the normal cone $N(x(t);C(t))$ is inactive, then
the first
component $\ou_{11}$ should employ the force as twice as the second one
$\ou_{12}$ does. In this case we get $\ou_{11}=-\frac{1}{3}$ and
$\ou_{12}=\frac{1}{6}$ with the minimum cost $\frac{53}{24}\approx 2.208$. It
seems that we might be on the right track of finding the optimal solution, but
this solution turns out to be nonoptimal if the normal vector $\eta_1(-1,-2)$ is
nontrivial. In the latter case the normal cone and the optimal control $\ou(t)$
provide the forces $\frac{1}{25t}(-1,-2)$ and $\(-\frac{2}{5},\frac{1}{10}\)$,
respectively. Thus the total force is
$\frac{1}{25}(-1,-2)+\(-\frac{2}{5},\frac{1}{10}\)=\(-\frac{9}{25},\frac{9}{50}
\)=\frac{9}{50}(-2,1)$, which actually points in the direction of vector
$(-2,1)$ and hence keeps the particle on the right track. More interestingly,
although the minimum cost is $\frac{441}{200}=2.205$ which is very close to the
cost in the former case, the contribution of the active normal vector
$\frac{1}{25}(-1,-2)$ still plays a very crucial role in reducing the
cost functional.
\end{remark}

\end{document}